\documentclass{siamltex}
\usepackage{amsfonts}
\usepackage{amssymb,amsmath}
\usepackage{graphicx}
\usepackage{color}
\usepackage{hyperref}

\DeclareMathAlphabet{\itbf}{OML}{cmm}{b}{it}

\def\EE{\mathbb{E}}

\def\eps{\varepsilon}

\def\bx{{\itbf x}}
\def\bw{{\itbf w}}

\def\om{{\omega}}
\def\la{{\lambda}}
\def\cD{{\mathfrak D}}
\def\cS{{\mathcal S}}
\def\cL{{\mathcal L}}

\def\cA{{\mathcal A}}
\def\cR{{\mathcal R}}
\def\cM{{\mathcal M}}
\def\cI{{\mathcal I}}

\def\cF{{\mathfrak F}}

\renewcommand{\hat}{\widehat}

\begin{document}

\title{A quantitative study of source imaging in random waveguides}

\author{
Liliana Borcea\footnotemark[1],
Josselin Garnier\footnotemark[2],
and
Chrysoula Tsogka\footnotemark[3]
}

\maketitle

\renewcommand{\thefootnote}{\fnsymbol{footnote}}

\footnotetext[1]{Computational and Applied Mathematics, Rice
  University, Houston, TX 77005. {\tt borcea@rice.edu}}
\footnotetext[2]{Laboratoire de Probabilit\'es et Mod\`eles
  Al\'eatoires \& Laboratoire Jacques-Louis Lions, Universit{\'e}
  Paris VII, 75205 Paris Cedex 13, France.  {\tt
    garnier@math.univ-paris-diderot.fr}} \footnotetext[3]{Applied
  Mathematics, University of Crete \& IACM/FORTH, 71409 Heraklion, Greece.  {\tt
    tsogka@tem.uoc.gr}}
\begin{abstract}
  We present a quantitative study of coherent array imaging of remote
  sources in randomly perturbed waveguides with bounded cross-section.
  We study how long range cumulative scattering by perturbations of
  the boundary and the medium impedes the imaging process. We show
  that boundary scattering effects can be mitigated with filters that
  enhance the coherent part of the data. The filters are obtained by
  optimizing a measure of quality of the image. The point is that
  there is an optimal trade-off between the robustness and resolution
  of images in such waveguides, which can be found adaptively, as the
  data are processed to form the image. Long range scattering by
  perturbations of the medium is harder to mitigate than scattering by randomly
  perturbed boundaries. Coherent imaging
  methods do not work and more complex incoherent methods, based on
  transport models of energy, should be used instead. Such methods are
  nor useful, nor needed in waveguides with perturbed boundaries. We
  explain all these facts using rigorous asymptotic stochastic
  analysis of the wave field in randomly perturbed waveguides.  We
  also analyze the adaptive coherent imaging method and obtain a
  quantitative agreement with the results of numerical simulations.
\end{abstract}
\renewcommand{\thefootnote}{\arabic{footnote}}

\section{Introduction}
\label{sect:intro}
\begin{figure}[t!]
    \vspace{-1.6in}
     \hspace{1.in}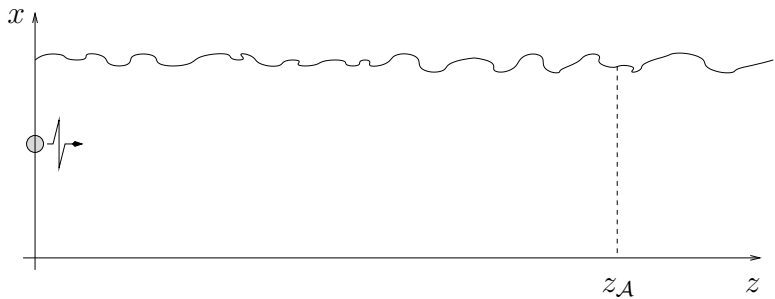
   \vspace{-0.05in}
   \caption{Schematic of the problem setup. A source emits a signal in
     a waveguide and the wave field is recorded at a remote array.
     The waves propagate along the range axis $z$. The waveguide is
     bounded in the cross-range direction $x$. The bottom boundary is
     rigid and flat. The pressure release top boundary may fluctuate.
     The system of coordinates has the origin of range at the source.
     The array is shown on the right of the source, at range $z_\cA$.}
   \label{fig:schem}
\end{figure}

We present a theoretical and numerical study of imaging remote sources
in random waveguides, using an array of sensors that record acoustic
waves. The waveguide effect is caused by the boundary of the
cross-section, which traps the waves and guides the energy along the
range direction $z$, as illustrated in Figure \ref{fig:schem}.  We
restrict our study to two-dimensional waveguides, because the
numerical simulations become prohibitively expensive in three
dimensions. The results are similar in three-dimensional waveguides
with bounded cross-section. We refer to \cite{BG-13} for an analysis
of wave propagation and imaging in three-dimensional random waveguides
with unbounded cross-section.

Scattering at the boundary creates multiple traveling paths of
the waves from the source to the receiver array.  Mathematically, we can write
the wave field $p$ (the acoustic pressure) as a superposition of a
countable set of waveguide modes, which are solutions of the
homogeneous wave equation.  Finitely many modes propagate in the range
direction at different speeds, and the remaining infinitely many modes
are evanescent waves that decay exponentially with range. We may
associate the propagating modes with planar waves that strike the
boundaries at different angles of incidence. The slow modes correspond
to near normal incidence. They reflect repeatedly at the boundary, thus
traveling a long path to the array. The fast modes correspond to
small grazing angles and shorter paths to the array.

In ideal waveguides with straight boundaries and wave speed that is
constant or varies smoothly with cross-range, the wave equation is
separable and the modes are uncoupled. In particular, each mode has a
constant amplitude which is determined by the source excitation. We
study perturbed waveguides with small and rapid fluctuations of the
boundaries and of the wave speed, due to numerous weak
inhomogeneities. Such fluctuations are not known and are of no
interest in imaging. However, they cannot be neglected because they
cause wave scattering that accumulates over long distances of
propagation. To address the uncertainty of the boundary and wave speed
fluctuations, we model them with random processes, and thus speak of
random waveguides. The array measures one realization of the random
field $p$, the solution of the wave equation in one realization of the
random waveguide. That is to say, for a particular perturbed boundary
and medium. When cumulative scattering by the perturbations is
significant, the measurements are quite different from those in ideal
waveguides. Furthermore, if we could repeat the experiment for many
realizations of the perturbations, we would see that the measurements
change unpredictably, they are statistically unstable.

The expectation (statistical mean) $\EE[p]$ of the wave is called the
\emph{coherent} field. This is the part of the data that is useful for
coherent imaging, because we can relate it to the unknown location of
the source, in spite of the uncertainty of the perturbations in the
waveguide. The challenge is to process the data in order to enhance
the coherent part $\EE[p]$ and mitigate the unwanted reverberations $p
- \EE[p]$, the \emph{incoherent} part.  Coherent methods without such
processing give images that are difficult to interpret and unreliable.
They change unpredictably with the realization of the random
waveguide, they are not statistically stable.

We refer to \cite{ABG-12} for a rigorous asymptotic stochastic
analysis of the wave field $p$ in waveguides with randomly perturbed
boundaries, and to \cite{kohler77,dozier,book07,garnier_papa} for
waveguides with randomly perturbed media. The analysis shows that $p$
can be modeled as a superposition of ideal waveguide modes that are
coupled by scattering at the random perturbations.  Explicitly, the
modes have amplitudes that are random functions of frequency and
range, and satisfy a coupled system of stochastic differential
equations. Their expectations decay exponentially with range, 
on mode- and frequency-dependent length scales called {\em scattering mean free
  paths}.  The decay means that the incoherent fluctuations of the
amplitudes gain strength, and once they become dominant, the modes
should not be used in coherent imaging.

It is not surprising that the scattering mean free paths are longer
for the fast propagating modes than the slower ones. This is because
the latter are waves that take longer trajectories from the source to
the array, and interact more with the perturbations of the boundaries
and the medium. We show in this paper that a successful imaging
strategy depends on which perturbations play the dominant role in the
waveguide.  If scattering from perturbed boundaries dominates, the
fast modes have a much longer scattering mean free path than the
slower modes.  Therefore, the data remain partially coherent at long
ranges and we can seek an adaptive imaging approach that detects the
slow modes with incoherent amplitudes and suppresses them.  The longer
the range, the fewer the modes that remain coherent, so there is a
trade-off between the statistical stability and the resolution of the
images, which can be optimized with the adaptive method.

When we compare the effect of perturbed boundaries to that of
perturbed media, for similar amplitude and correlation length of the
fluctuations, we find two essential differences: The latter gives much
shorter scattering mean free paths for the faster modes, and the rate
of change of these scales with the mode index is much slower.  There
is no trade-off between statistical stability and resolution of
coherent images in such waveguides. As the range increases, the mode
amplitudes become incoherent on roughly the same range scale, so there
is no gain in removing the slow modes. Coherent imaging fails and
should be replaced by incoherent methods, based on transport equations
for the energy resolved locally in time and over the modes i.e., over
the direction of propagation of the associated plane waves.  We refer
to \cite{BIT-10} for an example of incoherent imaging in random
waveguides. These methods are more complex and computationally
involved than the coherent ones. They are designed to work at ranges
that exceed the scattering mean free paths, but they also fail when
the source is further from the array than the \emph{equipartition
  distance}.  This is the range scale over which the energy of the
wave becomes distributed uniformly over the modes, independent of the
source excitation. The waves scatter so much while they travel this
distance that they lose all information of their initial state, thus
making imaging impossible. 

We show that in waveguides with interior inhomogeneities the
equipartition distance is much longer than the scattering mean free
path of the modes, so there is an observable range interval over which
coherent imaging fails, but incoherent imaging succeeds. This is not
the case for waveguides with perturbed boundaries where the
equipartition distance is almost the same as the scattering mean free
path of the fast modes. When coherent imaging fails in such
waveguides, no imaging method can succeed, so there is no advantage in
using the more complex, incoherent approaches.

The adaptive coherent imaging method proposed in this paper is based
on a figure of merit of the quality of the image, which accounts for
the trade-off between its statistical stability and resolution. There
are many such figures of merit. We choose one that is simple and
serves our purpose. In practice, it may be improved for example by
incorporating prior information about the support of the source
distribution. The method searches for weights of the data decomposed
over the waveguide modes, in order to optimize the figure of merit.
We apply the results of the asymptotic stochastic analysis in
\cite{kohler77,book07,garnier_papa,ABG-12} to derive theoretically the
weights, and show that they are in good agreement with those from the
numerical simulations in waveguides with random boundaries. We also
show that coherent imaging fails in random waveguides with interior
inhomogeneities, as predicted by the theory.

The paper is organized as follows. We begin in section \ref{sect:form}
with the formulation of the problem.  Then we describe in section
\ref{sect:model} the model of the array data in ideal and randomly
perturbed waveguides. The comparison of long range cumulative
scattering effects of boundary perturbations and interior
inhomogeneities is in section \ref{sect:scat}. The results motivate
the adaptive coherent imaging method described and analyzed in section
\ref{sect:imag}. The numerical simulations are in section
\ref{sect:num}.  We end with a summary in section \ref{sect:sum}.

We dedicate this work to George Papanicolaou on the occasion 
of his 70th birthday.

\section{Formulation of the source imaging problem}
\label{sect:form}
Consider a two-dimensional waveguide with range axis denoted by $z \in
\mathbb{R}$ and transverse coordinate (cross-range) $x$ belonging to a
bounded interval, the waveguide cross-section, as illustrated in
Figure \ref{fig:schem}. We assume a pressure release top boundary that
may be perturbed, and a flat and rigid bottom boundary.  Waveguides
with perturbations of both boundaries are studied in \cite{ABG-12}.

The pressure field $p(t,x,z)$ satisfies the wave equation
\begin{equation}
\label{we}
\left[\partial^2_{z} + \partial^2_x
  - \frac{1}{c^2(x,z)} \partial^2_t \right]p(t,x,z) = F(t,x,z) \, , 
\quad  x \in (0,\mathfrak{D}(z)), \quad z \in \mathbb{R}, \quad t>0,
\end{equation}
with boundary conditions 
\begin{equation}
  p(t,\mathfrak{D}(z),z) = \partial_x p(t,0,z) =0, \quad  z \in 
  \mathbb{R}, \quad t>0.
\label{eq:bc}
\end{equation}
Here $t$ is time, $c(x,z)$ is the wave speed, $\mathfrak{D}(z)$ is the
waveguide cross-section, and $F(t,x,z)$ models the source excitation.
In ideal waveguides the boundaries are straight
\begin{equation}
\mathfrak{D}(z) = D, \quad \forall z \in \mathbb{R},
\end{equation}
and the wave speed is independent of range. We take it equal to the
constant $c_o$. This simplification leads to explicit formulas in the
analysis of coherent imaging, but the results extend to speeds that
vary smoothly in $x$. In perturbed waveguides the boundary and the
wave speed have small amplitude fluctuations
\begin{equation}
|\mathfrak{D}(z) - D | \ll D \quad \mbox{and} \quad 
|c(x,z)-c_o| \ll c_o,
\label{eq:pert}
\end{equation}
modeled by random processes, as explained in sections \ref{sect:RB}
and \ref{sect:RM}.

We study the point spread function of coherent imaging methods,
so we let 
\begin{equation}
F(t,x,z) = e^{-i \om_o t} f(Bt) \delta(x - x_o ) \delta(z)  \,,
\label{eq:source}
\end{equation}
with the origin of the range axis at the point-like source, with
cross-range coordinate $x_o$.  The emitted signal is a pulse, modeled
by function $f$ of dimensionless arguments, with Fourier transform
$\hat f$ supported in the interval $[-\pi,\pi]$. The multiplication by
the carrier oscillatory signal $e^{-i \om_o t}$ centers the support of
the Fourier transform of the pulse at $\om_o$,
\begin{equation}
  \int_{-\infty}^\infty dt \, e^{-i \om_o t} f(Bt) e^{i \om t}  = 
  \frac{1}{B} \hat f\left(\frac{\om-\om_o}{B}\right).
\label{eq:signal}
\end{equation}
Therefore, the angular frequency $\om$, the dual variable to $t$,
belongs to the interval $[\om_o -\pi B, \om_o + \pi B]$, where
$\om_o/(2\pi)$ is the central frequency, and $B$ is the bandwidth.

The array is a collection of sensors that are far away from the
source, at range $z_\cA$, and record the pressure field $p$. The
recordings are the array data. The goal of coherent imaging is to
superpose the data after proper synchronization and weighting, in
order to form an imaging function. The synchronization is relative to
a search point that sweeps a search domain where we seek the source.
It amounts to solving backward the wave equation in the ideal
waveguide, with the source at the array and the emitted signal given
by the time reversed data. The imaging function is defined point-wise
by the resulting solution at the search point. This process is called
back-propagation.

A useful imaging function has the following qualities: (1) It peaks
near the unknown source.  (2) It is negligible away from the source.
The smaller the domain where it is large, the better the resolution.
(3) It is robust with respect to the unknown perturbations in the
waveguide.

Coherent imaging can succeed in random waveguides up to ranges where
the array data maintain some coherence. The asymptotic stochastic
theory developed in \cite{kohler77,book07,garnier_papa,ABG-12} allows
us to quantify the loss of coherence of the amplitudes of the
waveguide modes. We use the results to explain the limitations of
coherent imaging, and to motivate and analyze the adaptive imaging
approach.

\section{Model of the array data}
\label{sect:model}
We begin in section \ref{sect:homog} with the model of the data in
ideal waveguides. Then, we consider waveguides with a random pressure
release boundary in section \ref{sect:RB}, and with random wave speed
in section \ref{sect:RM}. The results extend to waveguides with both
types of random perturbations. We separate them in order to compare
their cumulative scattering effects on the imaging process.

\subsection{Ideal waveguides}
\label{sect:homog}
When the boundaries are flat and the wave speed is constant, the wave
equation is separable and we can write the solution as a superposition
of independent waveguide modes.  A waveguide mode is a monochromatic
wave $P(t,x,z) = \hat{P}(\omega,x,z) e^{- i \omega t}$, where
$\hat{P}(\omega,x,z)$ satisfies the Helmholtz equation
\begin{equation}
\label{eq:e0}
\left[{\partial_z^2} + {\partial_x^2} + k^2 \right]
\hat{P}( \omega,x,z) = 0 \, , \quad  x \in (0,D), ~ ~ z 
\in \mathbb{R},
\end{equation}
with boundary conditions
\begin{equation}
\label{eq:e1}
\hat P(\omega,D,z) = \partial_x \hat P(\omega,0,z) = 0, \quad z \in 
\mathbb{R},
\end{equation}
and radiation conditions as $|z| \to \infty$. Here $k = \om/c_o$ is
the wavenumber.

The linear operator $\partial^2_x + k^2$ defined on the vector space
of functions in $C^2(0,D)$ that vanish at $x=D$ and have zero
derivative at $x=0$, is self-adjoint in $L_2(0,D)$. Its spectrum
consists of a countable set of real and simple eigenvalues
$\{\lambda_j (\omega)\}_{j \ge 1}$, assumed sorted in descending
order. Because we assumed that $c_o$ is constant, we can write them
explicitly,
\begin{equation}
  \lambda_j(\om) = k^2 -  \left[\frac{\pi(j-{1}/{2})}{D}\right]^2  , 
\quad \quad j = 1, 2, \ldots.
\label{eq:eq2}
\end{equation}
The eigenfunctions form a complete orthonormal set in $L_2(0,D)$, and
are given by
\begin{equation}
\label{eq:e4}
\phi_j(x) = \sqrt{\frac{2}{D}} \cos \left[ \frac{\pi (j-{1}/{2}) x}{D}
\right], \qquad j = 1,2, \ldots.
\end{equation}

Note that only the first $N(\omega)$ eigenvalues are positive, where
\begin{equation}
N(\omega) = \left \lfloor
  {k D}/{\pi} +1/2 \right \rfloor,
\label{eq:e3}
\end{equation} 
and $\lfloor ~ \rfloor$ denotes the integer part.  They define the
modal wavenumbers $\beta_j(\omega)=\sqrt{\lambda_j(\omega)}$ of the
forward ($+$) and backward ($-$) propagating modes
\begin{equation}
\label{eq:e5}
\hat{P}_j^{(\pm)}(\omega,x,z) = \phi_j(x) e^{ \pm i
  \beta_j(\omega) z}, \quad j = 1, \ldots, N(\omega). 
\end{equation}
The remaining infinitely many modes are evanescent
\begin{equation}
  \hat{P}_j(\omega,x,z) = \phi_j(x) e^{ - 
    \beta_j(\omega) |z|}, \quad 
  j > N(\omega)\,, 
\label{eq:e6}
\end{equation}
with wavenumber $\beta_j (\omega)= \sqrt{-\lambda_j(\omega)}\, $.

\subsubsection{Plane wave analogy}
With the expression (\ref{eq:e4}) of the eigenfunctions, we can write
the forward propagating modes as
\begin{equation}
\label{eq:PW}
\hat P_j^{(+)}(\omega,x,z) = \frac{1}{\sqrt{2 D}} \left[ e^{i
    \left(\frac{\pi(j-1/2)}{D},\beta_j\right)\cdot (x,z)} + e^{i
    \left(-\frac{\pi(j-1/2)}{D},\beta_j\right)\cdot (x,z)} \right].
\end{equation}
A similar formula holds for the backward propagating modes, with a
negative sign in front of $\beta_j$. Equation (\ref{eq:PW}) shows that
the modes are associated with monochromatic plane waves that travel in
the direction of the slowness vectors
\[
{\itbf K}_j = \left( \pm \frac{\pi(j-1/2)}{D},\beta_j\right),
\]
and strike the boundaries where they reflect. The slowness vectors of
the first modes are almost parallel to the range axis,
\[ {\itbf K}_1 = \left( \pm \frac{\pi}{2 D}, \beta_1\right), \qquad
\frac{\pi}{2D} \approx \frac{k}{2 N} \ll \beta_1 \approx k,
\]
where the approximation is for a large $N(\om)$.  These waves travel
quickly to the array, at speed that is approximately equal to $c_o$.
The slowness vectors of the last modes are almost parallel to the $x$
axis
\[
{\itbf K}_N = \left( \pm \frac{\pi(N-1/2)}{D},\beta_N\right), 
\qquad \frac{\pi(N-1/2)}{D} \approx k \gg \beta_N.
\]
These waves strike the boundary many times, at almost normal
incidence. They propagate very slowly to the array, on a long 
trajectory.

\subsubsection{Data model}
\label{sect:homDM}
To simplify the analysis, we assume that the bandwidth is not too
large, so that
\begin{equation}
N(\om) = N(\om_o), \qquad \forall \om \in [\om_o-\pi B,\om_o+\pi B].
\label{eq:e7}
\end{equation}
We denote henceforth the number of propagating modes by $N$, without
any arguments. We also suppose that there are no standing waves in the
waveguide, which means that none of the wavenumbers $\beta_j$ vanish.

The pressure field for $z > 0$ is modeled by a superposition of
forward going and evanescent waves
\begin{equation*} 
  {p}(t,x,z) = \int \frac{d \om}{2 \pi} e^{-i \om
    t} \bigg[ \sum_{j=1}^N \frac{\hat{a}_{j,o}
    (\omega)}{\sqrt{\beta_j(\omega)}} e^{i\beta_j 
    (\omega)z} \phi_j(x)
  + \sum_{j=N+1}^\infty \hspace{-0.05in} \frac{\hat{e}_{j,o}
    (\omega)}{\sqrt{\beta_j(\omega)}} e^{- \beta_j(\omega) z}
  \phi_j(x) \bigg] \, .
\end{equation*}
The modes do not interact with each other, so their amplitudes are
independent of range. They are obtained from the source conditions 
\begin{align*}
\hat p(\om,x,0+) & = \hat p(\om,x,0-)\, , \\
\partial_z \hat p(\om,x,0+)-\partial_z \hat p(\om,x,0-)&= \frac{1}{B}\hat
f\left(\frac{\om-\om_o}{B}\right) \delta(x-x_o)\, ,
\end{align*}
which give
\begin{eqnarray}
  \hat{a}_{j,o}(\omega) &=&   \frac{ \phi_j( x_o )}{2i B
    \sqrt{\beta_j(\omega)}} \hat{f}\left(\frac{\omega-
      \omega_o}{B}\right)\,,
  \quad j =
  1, \ldots, N, \nonumber \\ \hat{e}_{j,o}(\omega) &=&
  - \frac{ \phi_j( x_o )}{2 B
    \sqrt{\beta_j(\omega)}} \hat{f}\left(\frac{\omega-
      \omega_o}{B}\right)\,, \quad j > N.
\label{eq:idealab} 
\end{eqnarray}

The model of the array data is given by 
\begin{equation} 
  {p}(t,x,z_\cA) \approx \sum_{j=1}^N \int \frac{d \om}{2 \pi B} 
  \hat{f}\left(\frac{\omega-
      \omega_o}{B}\right)\frac{\phi_j(x_o)}{2 
    i \beta_j(\om)}\, \phi_j(x)\, e^{i\beta_j 
    (\omega)z - i \om t}.
\label{eq:homop}
\end{equation}
The approximation is because we neglect the evanescent modes at 
the large range $z_\cA$ of the array.

\subsection{Waveguides with randomly perturbed boundary}
\label{sect:RB}
The pressure release boundary has small fluctuations around the value
$D$
\begin{equation}
  \cD(z) = D\left[1 +\nu\left(\frac{z}{\ell}\right)\right],
\label{eq:RB}
\end{equation}
where $\nu$ is a zero mean random process of dimensionless arguments.
We assume that it is stationary and mixing, which means in particular
that its covariance function
\begin{equation}
\label{eq:covar}
\cR_\nu(\zeta) = \EE\left[\nu(0) \nu(\zeta)\right]
\end{equation}
is integrable over the real line. The scaling of the argument of $\nu$
in (\ref{eq:RB}) indicates that the fluctuations are on the length
scale $\ell$, the correlation length. 

Let $\eps$ be the small parameter that scales the amplitude of the
fluctuations $\nu$, defined by
\begin{equation}
\cR_\nu(0) = \eps^2 \ll 1.
\label{eq:defEps}
\end{equation}
The asymptotic analysis in \cite{ABG-12} is with respect to $\eps$,
in the scaling regime
\begin{equation}
\label{eq:scaleRB1}
\ell \sim \la_o,
\end{equation}
where $\la_o$ is the reference, order one length scale.  In this
regime the waves interact efficiently with the random perturbations,
but because their amplitude is small, their cumulative scattering
effect is observable only at long ranges. It is shown in \cite{ABG-12}
that the scaling for studying the transition from coherent to
incoherent waves should be
\begin{equation}
\eps^2 z_\cA \sim \la_o. 
\label{eq:scaleRB2}
\end{equation}

We recall directly from \cite{ABG-12} the model of the pressure
field
\begin{equation}
{p}\left(t,x,z_\cA\right) \approx \int 
\frac{d \om}{2 \pi} \sum_{j=1}^N
\frac{\hat{a}_{j}  (\omega,z_\cA)}{\sqrt{\beta_j(\omega)}}\phi_j(x) \, 
e^{i\beta_j
  (\omega)z_\cA-i \om t}   \, .
\label{eq:randp}
\end{equation}
It is similar to equation (\ref{eq:homop}), except that the mode
amplitudes are random functions of frequency and range $z_\cA$.
They are analyzed in detail in \cite{ABG-12}. Here we need only 
their first and second moments: 

The mean mode amplitudes are 
\begin{equation}
  \EE [ \hat{a}_j(\omega,z_\cA)]\approx \frac{ \phi_j( x_o )}{2i B
    \sqrt{\beta_j(\omega)}} \hat{f}\left(\frac{\omega-
      \omega_o}{B}\right)
  \exp \left[ - \frac{z_\cA}{\cS_j(\omega)} + 
    i \frac{z_\cA}{\cL_j(\omega)} \right],
\label{eq:meanA}
\end{equation}
where the approximation indicates that there is a vanishing residual
in the limit $\eps \to 0$. We recognize the first factor in
(\ref{eq:meanA}) as $\hat a_{j,o}$, the $j-$th mode amplitude in ideal
waveguides. However, $\EE[ \hat{a}_j(\omega,z_\cA)]$ decays
exponentially with $z_\cA$, on the length scale $\cS_j(\om)$ called
the \emph{scattering mean free path} of the $j-$th mode. It is given
by
\begin{equation}
  \frac{1}{\cS_j(\omega)} = \frac{\pi^4 \ell \, (j-1/2)^2}{ D^4 
    \beta_j(\om)} 
  \sum_{l=1}^N
  \frac{(l-1/2)^2}{\beta_l(\om)}
  \hat \cR_\nu\left[(\beta_j(\om)-\beta_l(\om))\ell\right],
\label{eq:SCMFP_B}
\end{equation}
in terms of the power spectral density $\hat \cR_\nu$, the Fourier
transform of the covariance $\cR_\nu$. We know that $\hat \cR_\nu \ge
0$ by Bochner's theorem, so all the terms in the sum are nonnegative.

Aside from the exponential decay, the mean amplitudes also display a
net phase that increases with $z_\cA$ on the mode-dependent length
scales $\cL_j(\om)$. We recall\footnote{Note that there is a
  typo in \cite[Eq.~(4.20)]{ABG-12}: there is no minus sign in the
  definition of $\Gamma^{(s)}_{jj}(\omega)$} its expression from
\cite{ABG-12}
\begin{align*}
  \frac{1}{\cL_j(\omega)} =& \, \frac{\pi^4 \ell \, (j-1/2)^2}{ D^4
    \beta_j(\om)} \cR_\nu(0) \sum_{l=1}^N
  \frac{(l-1/2)^2}{\beta_l(\om)} \gamma
  \left[\beta_j(\om)-\beta_l(\om)\right]
  \\
  &+ \frac{\pi^2(j-1/2)^2}{D^2 \beta_j(\om)}\cR_\nu(0)\left\{ -
    \frac{3}{2} + \sum_{l \ne j, l = 1}^N \frac{\left[\beta_l(\om) +
        \beta_j(\om) \right] (l-1/2)^2}{\beta_l(\om)(j+l-1)(j-l)}
  \right\}  \\
  &+\frac{\cR_\nu''(0)(j-1/2)^2}{\ell^2\beta_j(\om)} \left\{
    \frac{\pi^2}{6 } + \sum_{l \ne j, l = 1}^N
    \frac{\left[\beta_j(\om)-\beta_l(\om)\right]
      (l-1/2)^2}{\beta_l(\om)(j+l-1)^2(j-l)^2} \right\} +
  \kappa_j^{(e)}(\om),
\label{eq:PHASE}
\end{align*}
where 
\[
  \gamma(\beta) = 2 \int_0^\infty du \, \sin(\beta \ell u) \cR_\nu(u),
\]
and $\kappa_j^{(e)}(\om)$ is due to the interaction of the evanescent
waves with the propagating ones. It is given by 
\begin{align*}
  \kappa_j^{(e)} (\omega)=& \frac{2 \pi^4 (j-1/2)^2}{D^4 \beta_j(\om)}
  \sum_{l=N+1}^\infty \left\{ \frac{\ell \,(l-1/2)^2 }{\beta_l(\om)}
    \int_0^\infty du \, e^{-\ell \beta_l(\om) u } \cR_\nu(u) \cos
    \left[\ell \beta_j(\om) u \right]\right. \nonumber \\
  &\left.-\frac{(l-1/2)^2}{\beta_j^2(\om) + \beta_l^2(\om)} \right\} -
  \frac{2 \, \cR_\nu''(0)\, (j-1/2)^2}{\ell^2
    \beta_j(\om)}\sum_{l=N+1}^\infty \frac{(l-1/2)^2}{
    (l-j)^2(l+j-1)^2},
\end{align*}     
where we used integration by parts to simplify the formulas derived in 
\cite{ABG-12}.

The mean square mode amplitudes are
\begin{equation}
  \EE \left[ |\hat{a}_j(\omega,z_\cA)|^2 \right]\approx \frac{1}{4 B^2} 
  \left|\hat f \left(\frac{\om-\om_o}{B}
    \right) \right|^2
  \sum_{l=1}^N   \frac{\phi_l^2(x_o)}{\beta_l(\om)} T_{jl}(\om,z_\cA)
  \, ,
\label{eq:2ndM}
\end{equation}
with $N \times N$ matrix 
\[
{\bf T}(\om,z_\cA) = e^{ \boldsymbol{\Gamma}^{(c)} (\om) z_\cA} \, ,
\] 
and symmetric $N \times N$ matrix  $\boldsymbol{\Gamma}^{({c})}(\om)$ defined by
\begin{align}
  \Gamma_{jl}^{({c})}(\omega) =& \frac{\pi^4 \ell \, (j-1/2)^2
    (l-1/2)^2}{D^4 \beta_j(\om) \beta_l(\om)} \hat \cR_\nu\left[\ell
    (\beta_j(\om)-\beta_l(\om))\right],
  \quad j \neq l,\nonumber \\
  \Gamma_{jj}^{({c})}(\omega) =&- \sum_{l\neq j,l=1}^N
  \Gamma_{jl}^{({c})} (\om) , \quad j =1,\ldots,N. \label{eq:Gamma}
\end{align}
Let $\Lambda_j(\om)$ be the eigenvalues of $\boldsymbol{\Gamma}^{(c)}$, in
descending order, and ${\itbf u}_j(\om)$ its orthonormal eigenvectors.
We have from the conservation of energy that
\[
\Lambda_j(\om) \le 0,
\] 
so the limit $z_\cA \to \infty$ of the matrix exponential
\[
{\bf T}(\om,z_\cA)  = e^{\boldsymbol{\Gamma}^{({c})}(\om) z_\cA} = \sum_{j=1}^N e^{\Lambda_j(\om) z_\cA}
{\itbf u}_j(\om) {\itbf u}_j^T(\om),
\]
is determined by the null space of $\boldsymbol{\Gamma}^{(c)}(\om)$. Under the
assumption that the power spectral density $\hat \cR_\nu$ does not
vanish for any of the arguments in (\ref{eq:Gamma}), $\boldsymbol{\Gamma}^{(c)} (\om)$ is
a Perron-Frobenius matrix with simple largest eigenvalue $\Lambda_1(\om) =
0$.  The leading eigenvector is given by
\[
{\itbf u}_1 = \frac{1}{\sqrt{N}}(1,\ldots, 1)^T,
\] 
and as $z_\cA$ grows, 
\begin{equation}
\label{eq:LR}
\sup_{j,l=1,\ldots,N}\left| T_{jl}(\om,z_\cA) -
  \frac{1}{N} \right| \le O\left(e^{-\Lambda_2(\om) z_\cA}\right).
\end{equation}
Thus, the right handside in (\ref{eq:2ndM}) converges to a constant
\begin{equation}
\label{eq:limEq}
  \sum_{l=1}^N  \frac{\phi_l^2(x_o)}{\beta_l(\om)} T_{jl}(\om,z_\cA) \stackrel{z_\cA 
    \to \infty}{\longrightarrow} \frac{1}{N}  \sum_{l=1}^N  
\frac{\phi_l^2(x_o)}{\beta_l(\om)},
\end{equation}
on the length scale 
\[
\cL_{\rm equip} = 
-1/\Lambda_2(\om),
\] 
called the \emph{equipartition distance}. It is the range scale over
which the energy becomes uniformly distributed over the modes,
independent of the source excitation.

Equations (\ref{eq:meanA}), (\ref{eq:2ndM}) and (\ref{eq:limEq}) give
that the SNR (signal to noise ratio) of the amplitude of the $j-$th
mode satisfies
\begin{equation}
 {\rm SNR}[\hat a_j(\om,z_\cA)] = \frac{ \left| \EE[ 
      \hat a_j(\om,z_\cA)] \right|}{
    \sqrt{\EE \left[ \big| \hat a_j(\om,z_\cA) - 
          \EE[ \hat a_j(\om,z_\cA)] \big|^2 \right]}} 
    \sim \exp \left[-\frac{z_\cA}{\cS_j(\om)} \right].
\end{equation}
Therefore, the $j-$th mode loses coherence on the range scale
$\cS_j(\om)$, the scattering mean free path. The scaling
(\ref{eq:defEps}) of the amplitude of the fluctuations $\nu$ implies
that
\[
\cS_j \sim \eps^{-2} \la_o, 
\]
so the loss of coherence can be observed at ranges of the order
$\eps^{-2} \la_o$, as stated in (\ref{eq:LR}).

\subsection{Waveguides  with random medium}
\label{sect:RM}
The boundaries in these waveguides are straight, but the wave speed
is perturbed as 
\begin{equation}
\label{eq:pertc}
\frac{1}{c^2(x,z)} = \frac{1}{c_o^2} \left[ 1 + \mu\left(
    \frac{x}{\ell},\frac{z}{\ell}\right) \right]. 
\end{equation}
Here $\mu(x,z)$ is a mean zero, statistically homogeneous
random process of dimensionless arguments, with integrable
autocorrelation
\begin{equation}
\cR_\mu(\xi,\zeta) = \EE\left[\mu(0,0) \mu(\xi,\zeta)\right].
\label{eq:Rmu}
\end{equation}
As in the previous section, we model the small amplitude of the
fluctuations using the small dimensionless parameter $\eps$ defined by
\begin{equation}
\cR_\mu(0,0) = \eps^2 \ll 1.
\end{equation}
The scaling by the correlation length $\ell$ of both arguments of
$\mu$ indicates that the fluctuations are isotropic. We assume like
before that $ \ell \sim \la_o,$ and use the same long range scaling
(\ref{eq:scaleRB2}) to study the loss of coherence of the waves due to
cumulative scattering in the random medium.

The model of the array data,
the mean and intensity of the mode amplitudes look the same as
(\ref{eq:randp}), (\ref{eq:meanA}) and (\ref{eq:2ndM}), but the
scattering mean free paths $\cS_j(\om)$, the net phases
$\cL_j(\omega)$ and the matrix $\boldsymbol{\Gamma}^{({c})}(\omega)$ are different.  We
recall their expression from \cite[Chapter 20]{book07}.

The scattering mean free path of the $j-$th mode is given by
\begin{equation}
\label{eq:scRm}
\frac{1}{\cS_j(\om)} = \frac{k^4 \ell }{8 \beta_j(\om)} \sum_{l=1}^N 
\frac{1}{\beta_l(\om)} \hat \cR_{\mu_{jl}}\left[\left( 
\beta_j(\om)-\beta_l(\om)\right)\ell \right],
\end{equation}
where $\hat \cR_{\mu_{jl}}$ is the power spectral density of the
stationary process 
\begin{equation}
  \mu_{jl}(\zeta) = \int_0^D d x \, \phi_j(x) \phi_l(x) \mu 
  \left(\frac{x}{\ell},\zeta\right),
\end{equation}
with autocorrelation 
\begin{equation}
\cR_{\mu_{jl}}(\zeta) = \EE\left[\mu_{jl}(0) \mu_{jl}(\zeta)\right].
\end{equation}

The net phase of the $j-$th mode is  
\begin{equation}
\frac{1}{\cL_j(\om)} = \frac{k^4 \ell }{8 \beta_j(\om) } 
\sum_{l=1}^N \frac{1}{\beta_l(\om)} \gamma_{jl} \left[ 
\beta_j(\om)-\beta_l(\om) \right] + \kappa_j^{(e)}(\om),
\label{eq:phRm}
\end{equation}
where 
\begin{equation}
\gamma_{jl}(\beta) = 2 \int_0^\infty d u \, \sin ( \beta \ell u) 
\cR_{\mu_{jl}}(u),
\end{equation}
and the last term is due to the interaction of the evanescent modes
with the propagating ones
\begin{equation}
  \kappa_{jl}^{(e)}(\om) = \frac{k^4 \ell}{2 \beta_j(\om)} 
  \sum_{l = N+1}^\infty \frac{1}{\beta_l(\om)} \int_0^\infty 
du \,  e^{-\beta_l(\om) u} 
  \cR_{\mu_{jl}}(u) \cos \left[ \ell \beta_j(\om) u\right].
\end{equation}

The matrix $\boldsymbol{\Gamma}^{({c})}(\om)$ is symmetric, with entries given by 
\begin{align}
  \Gamma_{jl}^{({c})}(\omega) &= \frac{k^4 \ell
  }{8\beta_j(\om)\beta_l(\om)} \hat \cR_{\mu_{jl}}
  \left[ \left(\beta_j(\om)-\beta_l(\om)\right)\ell\right],  
\quad j \neq l, \nonumber \\
  \Gamma_{jj}^{({c})}(\omega) &=- \sum_{l\neq j,l=1}^N
  \Gamma_{jl}^{({c})}(\om) , \quad j=1,\ldots,N. \label{eq:GammaRM}
\end{align}
As before, we denote its eigenvalues by $\Lambda_j(\om) \le 0$, and
its orthonormal eigenvectors by ${\itbf u}_j$, for $j = 1, \ldots, N$.
Moreover, assuming that the power spectral density $\hat
\cR_{\mu_{jl}}$ does not vanish at any of the arguments
$(\beta_j-\beta_l)\ell$, we obtain from the Perron-Frobenius theorem
that the null space of $\boldsymbol{\Gamma}^{({c})}(\om)$ is one-dimensional and spanned
by
\[
{\itbf u}_1 = \frac{1}{\sqrt{N}}(1,1,\ldots, 1)^T.
\]
The long range limit of the matrix exponential is as in (\ref{eq:LR}),
and the equipartition distance is given by $ -1/\Lambda_2(\om)$, in
terms of the largest non-zero eigenvalue of $\boldsymbol{\Gamma}^{({c})}(\om)$.

\section{Comparisson of cumulative scattering effects}
\label{sect:scat}

It is not difficult to see by inspection of formulas
(\ref{eq:SCMFP_B}) and (\ref{eq:scRm}) that the scattering mean free
paths $\cS_j$ and the net phase range scales $\cL_j$ decrease
monotonically with the mode index. To obtain a quantitative
comparison of the net scattering effects of boundary and medium
perturbations, we consider here and in the numerical simulations two
examples of autocorrelations of the fluctuations $\nu(\zeta)$ and
$\mu(\xi,\zeta)$.  The conclusions drawn below extend qualitatively to
all fluctuations, but obviously, the scales depend on the expressions
of $\cR_\nu$ and $\cR_\mu$, the depth of the waveguide and the
correlation length relative to $\la_o$.

\begin{figure}[t]
\begin{center}
\begin{tabular}{c}
\includegraphics[width=8.0cm]{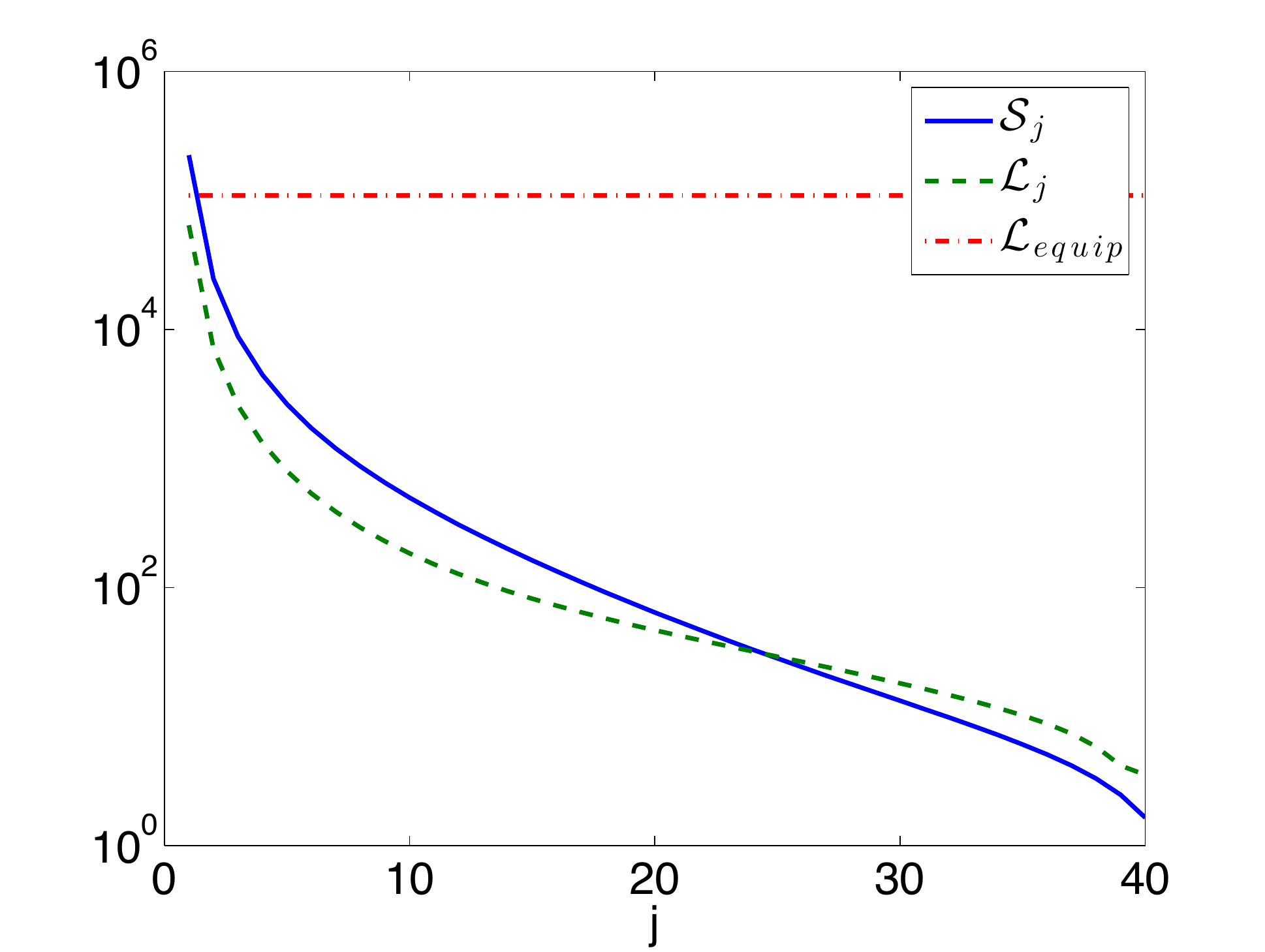} 
\end{tabular}
\end{center}
\caption{The characteristic  scales for a waveguide with
  random boundary.  Here $D=20 \lambda_o$, $\ell=\lambda_o/\sqrt{5}$
  and $\eps=0.013$. The abscissa is mode index and the ordinate is in
  units of $\la_o$. }
\label{fig:1}
\end{figure}

We take henceforth $D = 20 \la_o$, so that $N = 40$.  The
autocorrelation of the boundary fluctuations is of the so-called
Mat\'{e}rn$-7/2$ form
\begin{equation}
  \cR_\nu(\zeta)  = {\eps^2} 
  \left( 1 + |\zeta| +\frac{6 \zeta^2}{15} +\frac{
      |\zeta|^3}{15} \right) e^{-|\zeta|},
\label{eq:Rnu}
\end{equation}
with power spectral density 
\begin{equation}
  \hat \cR_\nu(\beta \ell) = \frac{32 \eps^2 }{5 \left[1 + 
(\beta \ell)^2\right]^4}.
\end{equation}
The correlation length is $\ell = \la_o/\sqrt{5}$, and the amplitude
of the fluctuations is scaled by $\eps = 0.013$. The
characteristic scales $\cS_j$, $\cL_j$ and the equipartition distance
$\cL_{\rm equip}$ are plotted in Figure \ref{fig:1}.
 
\begin{figure}[t]
\begin{center}
\begin{tabular}{c}
\includegraphics[width=8.0cm]{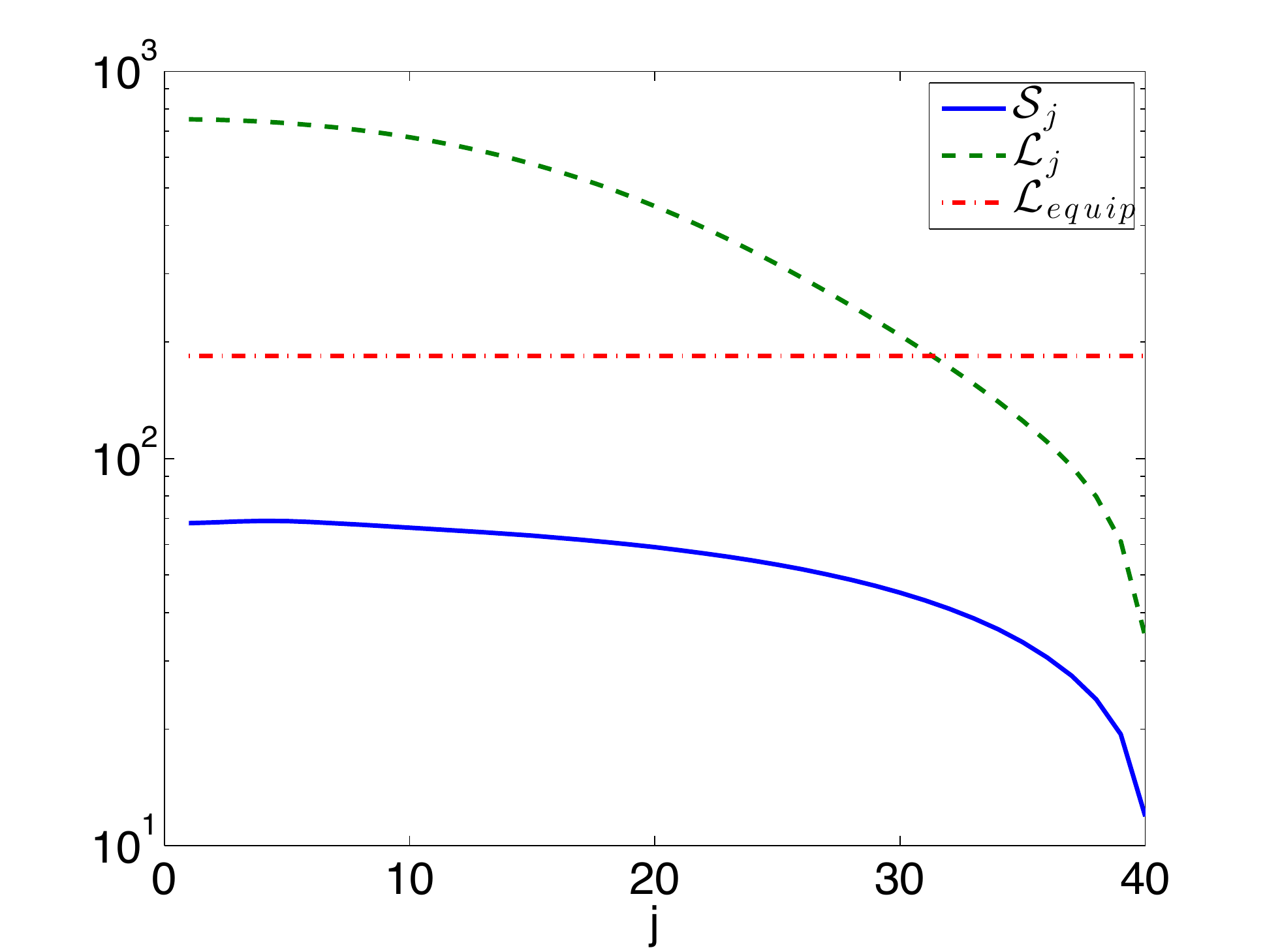} 
\end{tabular}
\end{center}
\caption{The characteristic scales for a waveguide with
  random medium. Here $D=20 \lambda_o$, $\ell=\lambda_o$, and
  $\eps=0.04$.  The abscissa is mode index and the ordinate is in
  units of $\la_o$.  }
\label{fig:2}
\end{figure}

The medium fluctuations have the Gaussian autocorrelation 
\begin{equation}
\cR_\mu(\xi,\zeta) = \eps^2 e^{-\frac{\xi^2 + \zeta^2}{2}},
\label{eq:RmuG}
\end{equation}
with correlation length $\ell = \la_o$ and amplitude scaled by
$\eps = 0.04$. The characteristic scales $\cS_j$,
$\cL_j$ and the equipartition distance $\cL_{\rm equip}$ are plotted
in Figure \ref{fig:2}.

We see in Figure \ref{fig:1} that the fast modes (with small index)
have much larger scattering mean free paths than the slow ones in the waveguides with
perturbed boundaries. When the arrays is at range $z_\cA \sim 100
\la_o$, roughly half of the mode amplitudes remain coherent, and we
can expect imaging to succeed if we filter out the slower modes, with
index $j > 20$. As $z_\cA$ increases, fewer and fewer modes remain
coherent, and imaging should become more difficult. Once $z_\cA$
exceeds the equipartition distance, which is similar to $\cS_1$ in
Figure \ref{fig:1}, imaging becomes impossible, because the wave field
forgets all the information about its initial state. Thus, when the
wave field loses all its coherence, no imaging method can succeed in
these waveguides.

Figure \ref{fig:2} shows that in media with random perturbations the
scattering mean free paths of the fast modes are shorter, and that
they decrease at a much slower rate with the mode index.  No mode
filtering can make coherent imaging succeed for $z_\cA \gtrsim 50
\la_o \sim \cS_1$, because all the mode amplitudes are incoherent.
Since the equipartition distance is much larger than $\cS_1$,
incoherent imaging is useful in these waveguides, in the range
interval
\[
50 \la_o \sim \cS_1 \lesssim z_\cA \le \cL_{\rm equip}
\sim 200 \la_o.
\]

\section{Adaptive coherent imaging}
\label{sect:imag}
We begin in section \ref{sect:DefIM} with the formulation of the
adaptive coherent imaging function. It models the backpropagation of
the weighted time reversed data to search points in a fictitious ideal
waveguide. The weights are chosen by optimizing a figure of merit of
the image. We calculate them explicitly in ideal and random
waveguides, in sections \ref{sect:ideal} and \ref{sect:random},
respectively.

\subsection{Coherent imaging}
\label{sect:DefIM}
The data are collected at the array, with sensors located in the set
$\{ \bx_r=(x_r,z_\cA), \quad r=1,\ldots,N_R\}$. The standard coherent
imaging function is given by 
\begin{equation} 
  \label{eq:IM1}
  {\mathcal I}(\bx) = \int_{-\infty}^{\infty} \frac{d \om}{2 \pi}
  \,  \frac{1}{N_R} \sum_{r=1}^{N_R} \overline{ \hat{p}(\omega,\bx_r) }
  \, \hat{G}_o(\omega,\bx_r,\bx),
\end{equation}
where $\bx=(x,z)$ are points in a search domain containing the unknown
location $\bx_o = (x_o,0)$ of the source, and $\hat{G}_o$ is the
outgoing Green's function in the ideal waveguide. It models the
propagation from $\bx_r$ to $\bx$, of the time reversed array data
with Fourier transform $\overline{ \hat p(\om,\bx_r)}$, where the bar
denotes complex conjugation.

We see from (\ref{eq:homop}) that 
\begin{equation}
\label{eq:IM2}
\hat{G}_o(\omega,\bx_r,\bx)
 = \sum_{j=1}^{N} \frac{\phi_j(x_r)}{2 i\beta_j(\omega)}\phi_j(x)
 e^{i \beta_j(\omega) (z_\cA-z)},
\end{equation}
so we can rewrite (\ref{eq:IM1}) as 
\begin{equation} 
  \label{eq:IM3}
  {\mathcal I}(\bx) = \int_{-\infty}^\infty \frac{d \om}{2 \pi} \,  
\sum_{j=1}^{N}
  \frac{1}{2 i \beta_j(\omega)} \overline{ \hat{p}_j(\omega,z_\cA) }
  \phi_j(x) e^{i \beta_j(\omega) (z_\cA-z)}
\end{equation}
with
\begin{equation}
\label{eq:IM4}
\hat{p}_{j}(\omega,z_\cA) = \frac{1}{N_R}  
\sum_{r=1}^{N_R}  \hat{p}(\omega,\bx_r)  \phi_j(x_r).
\end{equation}

The adaptive coherent imaging function is a modification of
(\ref{eq:IM3})
\begin{equation}
  \label{eq:IM5}
  {\mathcal I}(\bx;\bw) = \int_{-\infty}^\infty \frac{d \om}{2 \pi} \, 
  \sum_{j=1}^{N}
  \frac{w_j}{2 i \beta_j(\omega)} \overline{\hat{p}_j(\omega,z_\cA) }
  \phi_j(x) e^{i \beta_j(\omega) (z_\cA-z)},
\end{equation}
with data components $\hat p_j(\omega,z_\cA)$ weighted
by the entries in the complex vector  
\[
\bw = (w_1, \ldots, w_N)^T \in \mathbb{C}^N,
\]
with Euclidian norm 
\begin{equation}
\|\bw\| = \sqrt{\sum_{j=1}^N |w_j|^2} = 1.
\label{eq:norm}
\end{equation}

\subsubsection{Weight optimization}
\label{sect:WOPT}
We wish to optimize the weights so as to maximize the ratio of the
peak amplitude of the image normalized by its $L^2$-norm,
\begin{equation}
\label{eq:FIG_MERRIT}
\bw^\star = {\rm argmax}_{\bw \in \mathbb{W}} \cM(\bw), \quad \qquad
\cM(\bw) = \frac{\left| \cI(\bx^\star;\bw)\right|^2}{\|\cI(\cdot;\bw)
\|^2},
\end{equation}
where 
\[
  \mathbb{W} = \Big\{ \bw = (w_1, \ldots, w_N)^T \in \mathbb{C}^N, 
\quad \sum_{j=1}^N |w_j|^2 = 1 \Big\},
\]
and
\[
\|\cI(\cdot;\bw) \|^2 = \int_0^D d x \int_{-\infty}^\infty dz \, 
\left| \cI(\bx;\bw) \right|^2.
\]
The peak location $\bx^\star$ is expected to be at $\bx_o$, where the
source lies, and the optimization intends to focus the image around
it.  This is certainly true in ideal waveguides. In random waveguides
we need to ensure that the image is robust with respect to the unknown
perturbations. If this is not so, the image will have spurious peaks.

There are two requirements for obtaining robust images: The first is
that only the modes that are coherent contribute to the image.  Thus,
the weights should null the modes with scattering mean free paths that
are shorter than the range of the array. The second is that the
bandwidth be much larger than the decoherence frequency of the data.
This ensures that the incoherent part of the data averages out when we
integrate over the frequencies, like in the law of large numbers.

It is shown in \cite{ABG-12} and \cite[Chapter 20]{book07} that in our
regime the decoherence frequency is very small, of the order $\eps^2
\om_o$. Therefore, it is possible to have a bandwidth that is small
with respect to the central frequency, as assumed in section
\ref{sect:homDM}, and large with respect to the decoherence frequency.

As long as the two requirements above hold, we can analyze the optimal
weights using the theoretical figure of merit
\begin{equation}
  \cM_{\rm th}(\bw) = \frac{\left| \EE 
      \left[\cI(\bx_o;\bw)\right]\right|^2}{\EE \left[\|\cI(\cdot;\bw)
      \|^2\right]}.
\label{eq:ThM}
\end{equation}

\subsubsection{Simplifying assumptions}

In the analysis we suppose that the recordings of the acoustic
pressure are over an infinitely long time window, and approximate the
array by a continuum aperture, so that in the imaging function we can
replace sums over the sensors by integrals over the aperture. In
particular, we have 
\begin{equation}
\hat{p}_{j}(\omega,z_\cA) = \frac{1}{N_R}  
\sum_{r=1}^{N_R}  \hat{p}(\omega,\bx_r)  \phi_j(x_r) \approx 
\int_0^D dx \, 1_\cA(x) \hat p(\om,x,z_\cA) \phi_j(x),
\label{eq:ModeDec}
\end{equation}
where $1_\cA$ is the indicator function of the array. It is equal to
one in the cross-range support of the array and zero otherwise.  The
continuum approximation is valid when the sensors are close together,
at less than half a central wavelength $\la_o = 2 \pi c_o/\om_o$
apart.

We consider a full aperture array, spanning the entire cross-section
of the waveguide, so the indicator function $1_\cA$ in
(\ref{eq:ModeDec}) is identically one. The results extend to partial
apertures, but the formulas are more complicated and the optimal
weights are not easy to interpret. 

All these assumptions allow us to simplify the expression of the
imaging function, so that we can focus attention on the cumulative
scattering effects due to the random perturbations of the waveguide.

\subsection{Coherent imaging in ideal waveguides}
\label{sect:ideal}%
In this section we address the case in which the waveguide is ideal, i.e. without any random perturbation.

\subsubsection{Determination of the optimal weights}
We obtain from the model (\ref{eq:homop}) of the array data and 
the orthogonality of the eigenfunctions that 
\begin{equation}
\label{eq:IMh1}
\hat{p}_j(\omega) =\frac{1}{B}\hat f\left(\frac{\om-\om_o}{B}\right) 
\frac{\phi_j(x_o)}{2 i
  \beta_j(\omega)} e^{i \beta_j(\omega)z_\cA},
\end{equation}
and therefore
\begin{align}
  \cI(\bx;\bw) &= \frac{1}{4} \int_{-\infty}^\infty \frac{d \om}{2
    \pi B} \overline{\hat f\left(\frac{\om-\om_o}{B}\right)}
  \sum_{j=1}^{N} \frac{w_j}{\beta_j^2(\omega)}
  \phi_j(x)\phi_j(x_o) e^{- i \beta_j(\omega) z} \nonumber \\
  &\approx \frac{1}{4} \sum_{j=1}^{N} \frac{w_j}{\beta_j^2(\om_o)}
  \phi_j(x)\phi_j(x_o) \cF_j(z) \, .
\label{eq:IMh2}
\end{align}
Here we used that $B \ll \om_o$ and introduced the mode pulses
\begin{align}
\cF_j(z) &= \int_{-\infty}^\infty \frac{d \om}{2 \pi B}
  \overline{\hat f\left(\frac{\om-\om_o}{B}\right)} e^{- i
    \beta_j(\omega) z} \nonumber \\
&= \int_{-\infty}^\infty \frac{d h}{2 \pi}
  \overline{\hat f(h)} e^{- i
    \beta_j(\om_o + Bh) z} \nonumber \\
& \approx e^{-i \beta_j(\om_o) z} f\left[ - \beta'_j(\om_o) B
  z\right]
\label{eq:IMh3}
\end{align}
that peak at the range $z = 0$ of the source, with mode- 
and bandwidth-dependent resolution.  The modes propagate at speed
\begin{equation}
\label{eq:ModeSpeed}
\frac{1}{\beta'_j(\om_o)} = c_o \frac{\beta_j(\om_o)}{k_o}\, ,
\end{equation}
where $k_o = 2 \pi/\la_o$, and the range resolution of $\cF_j(z)$ is
determined by the distance traveled at this speed over the duration
$\sim 1/B$ of the pulse.

The focusing in cross-range is due to the summation over the modes.
Explicitly, when we evaluate (\ref{eq:IMh2}) at the range of the
source, we obtain
\begin{equation}
\cI((x,0);\bw) 
\approx \frac{f(0)}{4} \sum_{j=1}^{N} \frac{w_j}{\beta_j^2(\om_o)}
  \phi_j(x)\phi_j(x_o)\,.
\label{eq:IMh5}
\end{equation}
This is a sum of oscillatory terms unless $x = x_o$, so the image 
peaks at $x = x_o$, with resolution depending on the weights.

The figure of merit (\ref{eq:FIG_MERRIT}) is the ratio of the peak 
intensity 
\begin{equation}
  \left|\cI(\bx_o;\bw)\right|^2 \approx \frac{|f(0)|^2}{16} 
\Big|\sum_{j=1}^{N} \frac{w_j}{\beta_j^2(\om_o)}
  \phi_j^2(x_o)\Big|^2
\label{eq:IMh6}
\end{equation}
and the $L_2$ norm 
\begin{align}
  \nonumber \| \cI(\cdot ; \bw ) \|^2 &= \int_0^D dx
  \int_{-\infty}^\infty dz \left| \cI(\bx ; \bw )\right|^2 \nonumber
  \\&\approx \sum_{j=1}^{N} \frac{|w_j|^2\phi_j^2(x_o)}{16
    \beta_j^4(\om_o)}\hspace{-0.03in} \int_{-\infty}^\infty
  \hspace{-0.03in} \frac{dh}{2 \pi} \overline{\hat{f}(h)}
  \hspace{-0.03in}\int_{-\infty}^\infty \hspace{-0.03in} \frac{dh'}{2
    \pi} \hat{f}(h')\hspace{-0.03in}\int_{-\infty}^\infty
  \hspace{-0.1in}dz \, e^{ i [ \beta_j(\om_o + h B) -\beta_j(\om_o + B
    h') ]z}
  \nonumber \\
  &= \frac{\|f\|^2 c_o^2}{16 B \om_o} \sum_{j=1}^{N}
  \frac{|w_j|^2}{\beta_j^3(\om_o)} \phi_j^2(x_o) \, .
\end{align}
Here we used the orthonormality of the eigenfunctions, and relation
(\ref{eq:ModeSpeed}). We also introduced the notation
\[
\|f\|^2 = \int_{-\infty}^\infty \frac{d \omega}{2 \pi} \big|\hat
  f(\om)\big|^2 = \int_{-\infty}^\infty dt |f(t)|^2 .
\]
The figure of merit becomes 
\begin{equation}
  \cM(\bw) = C \frac{ \Big|\displaystyle \sum_{j=1}^{N} 
      \frac{w_j}{\beta_j^2(\om_o)}
      \phi_j^2(x_o)\Big|^2}{\displaystyle\sum_{j=1}^{N}
    \frac{|w_j|^2}{\beta_j^3(\om_o)} \phi_j^2(x_o)},
\label{eq:IMh7}
\end{equation}
with constant 
\[
C = \frac{B \om_o |f(0)|^2}{c_o^2 \|f\|^2}
\]
that plays no role in the optimization. Because $\cM$ is homogeneous
of degree zero in $\bw$, we can maximize $\cM(\bw)$ to obtain the
optimal $\bw^\star$ up to a multiplicative constant that we can then
determine from the normalization condition $\|\bw^\star\| = 1.$
The result is 
\begin{align}
  w^\star_j &= \frac{\beta_j(\om_o)}{\|{\boldmath{\beta}}\|_{x_o}},
  \quad j \in \mathbb{J}_{x_o} = \left\{ j = 1, \ldots, N, ~ ~ {\rm
      s.t.} ~ \phi_j(x_o) \ne 0 \right\},
  \label{eq:IMh8} \\
  w^\star_j &= 0, \quad j \in \mathbb{J}^c_{x_o} =\left\{1,
    \ldots, N\right\} \setminus \mathbb{J}_{x_o} , \nonumber
\end{align}
where we introduced the notation 
\[
\|{\boldmath{\beta}}\|_{x_o} := 
\sqrt{\sum_{j\in \mathbb{J}_{x_o}} \beta_j^2(\om_o)}\,.
\]
When the set $\mathbb{J}^c_{x_o}$ is empty, there is a unique
maximizer $\bw^\star$. Otherwise, there are infinitely many
maximizers, with arbitrary weights for mode indexes $j \in
\mathbb{J}^c_{x_o}$.  Equations (\ref{eq:IMh8}) define just one
solution.  Note however that all maxima of $\cM(\bw)$ are global
maxima, because the weights indexed by $j \in \mathbb{J}_{x_o}^c$
multiply $\phi_j(x_o) = 0$ in the figure of merit, and they play no
role in the behavior of the imaging function, given by 
\begin{equation}
  \cI(\bx;\bw^\star) \approx \frac{1}{4 \|\boldmath{\beta}\|_{x_o}}\sum_{j\in
    \mathbb{J}_{x_o}} 
  \frac{\phi_j(x)\phi_j(x_o)}{\beta_j(\om_o)} e^{-i \beta_j(\om_o)z}
  f\left[-\beta_j'(\om_o) B z\right].
\label{eq:IMh9}
\end{equation}

\begin{figure}[t]
\begin{center}
\begin{tabular}{c}
\includegraphics[width=4.2cm]{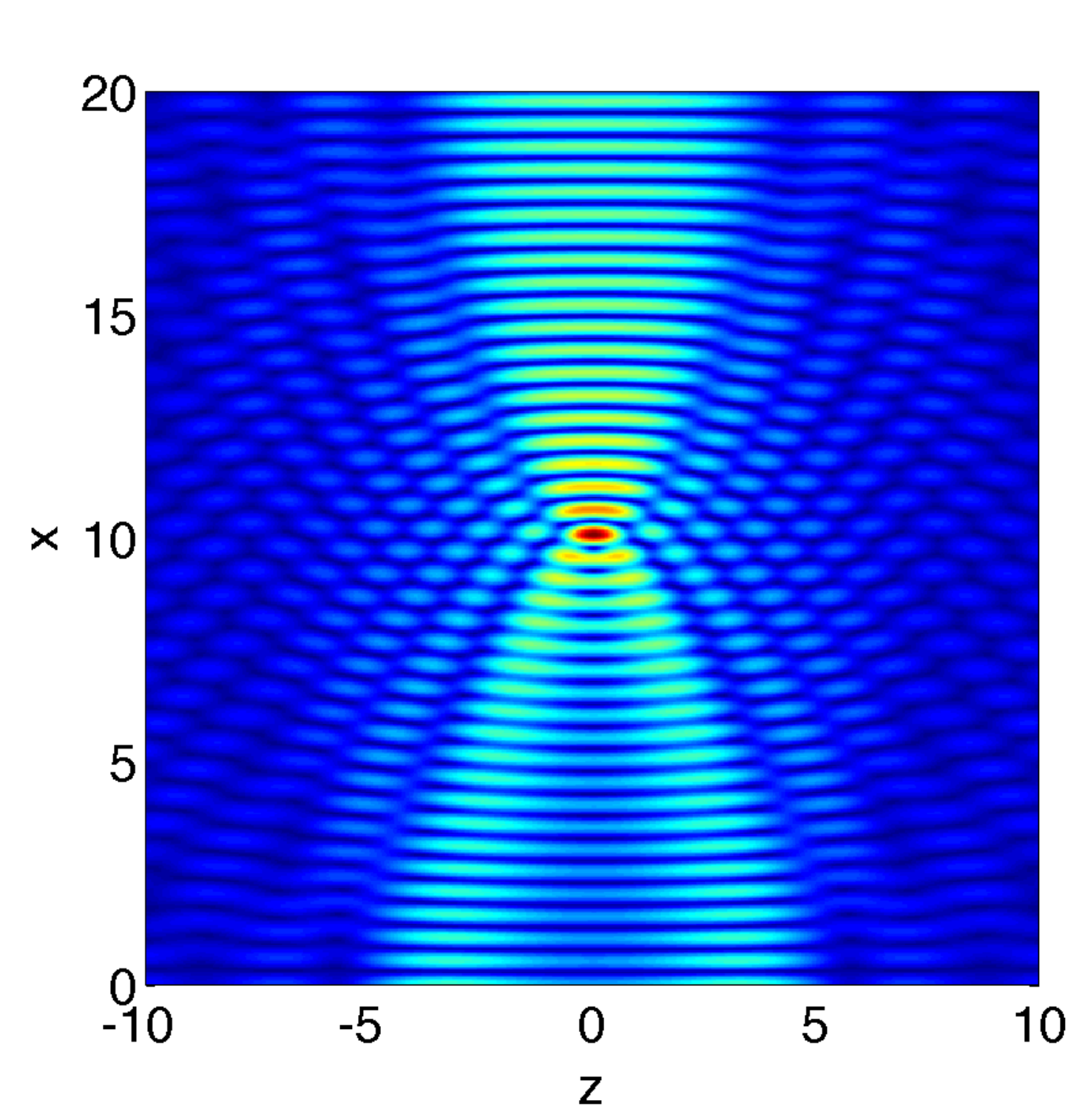}
\hspace{-0.1in}\includegraphics[width=4.2cm]{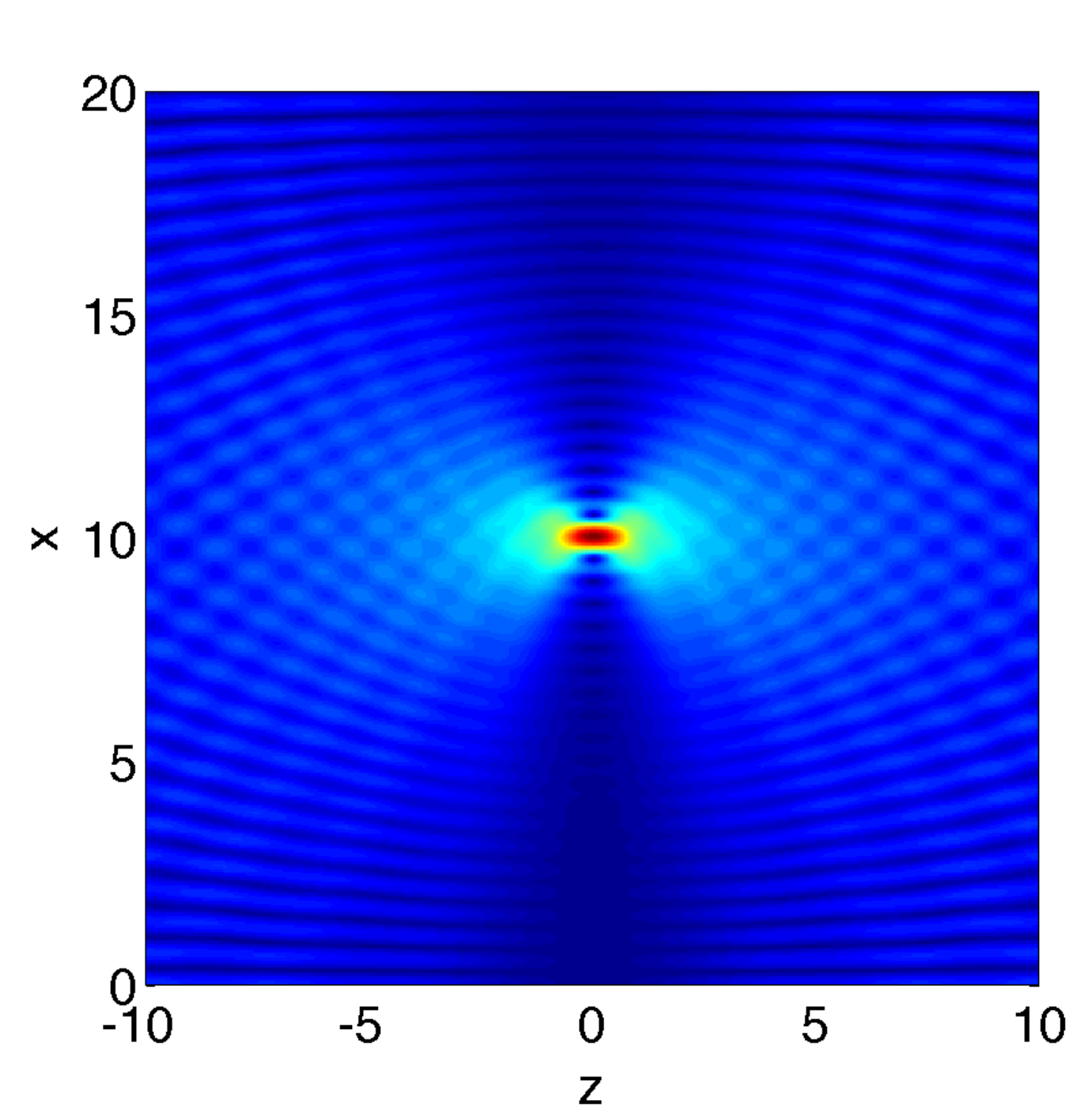}
\hspace{-0.1in}\includegraphics[width=4.2cm]{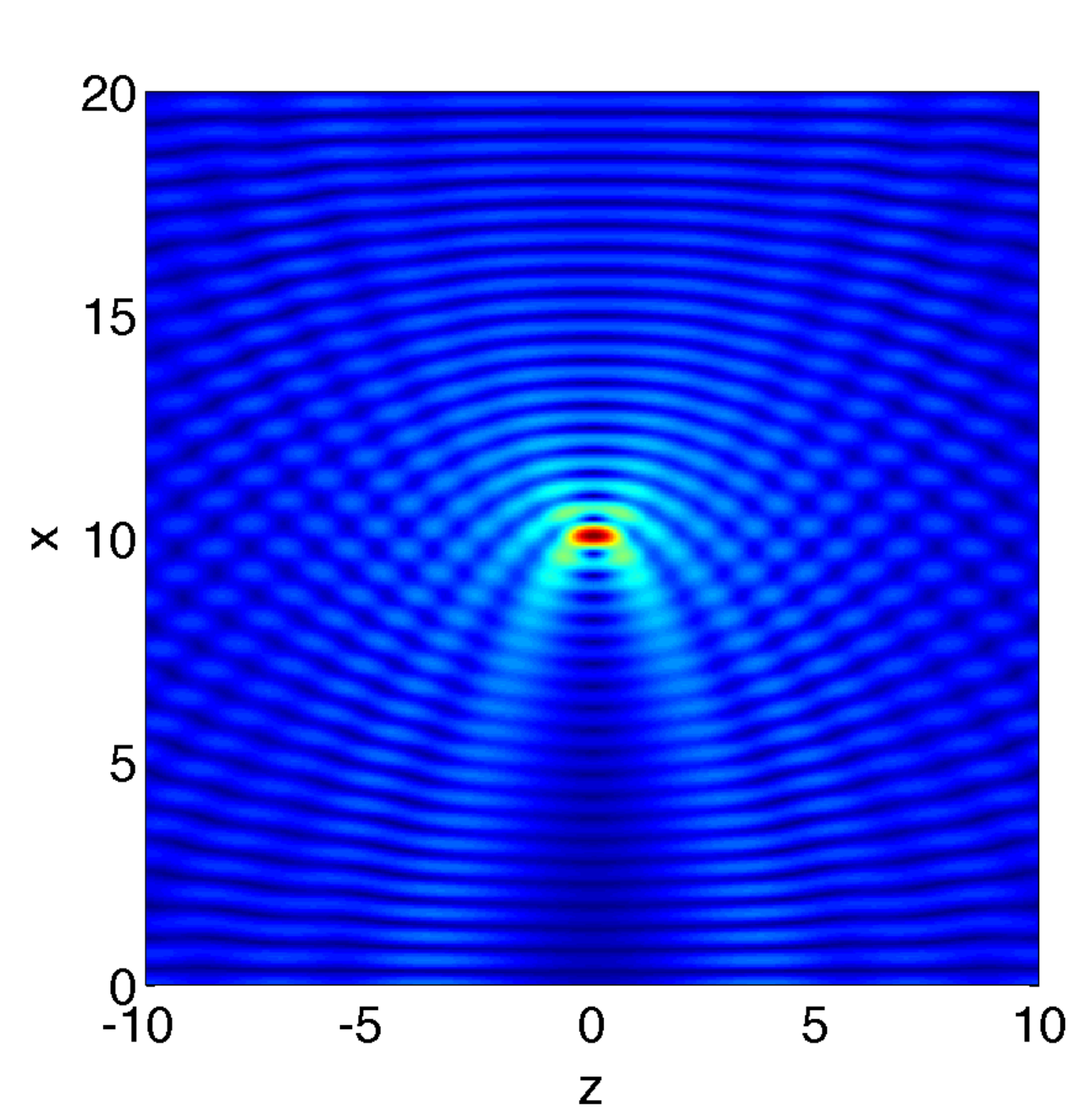}
\end{tabular}
\end{center}
\caption{The image $\cI(\bx;\bw)$ for weights $w_j = 1/\sqrt{N}$
  (left), $\bw = \bw^{\rm cr}$ (middle) and
  $\bw = \bw^\star$ (right). The abscissa is
  the range $z$ in $\lambda_o$ and the ordinate the cross-range in
  $\la_o$. Here $D = 20 \la_o$, and the pulse is defined in
  (\ref{eq:pulseG}).}
\label{fig:CR}
\end{figure}

\subsubsection{Discussion}
To motivate the figure of merit (\ref{eq:IMh7}) and illustrate the
effect of the optimization on the image, let us set $x_o = D/2$ and
consider a Gaussian pulse
\begin{equation}
f(u) = e^{-u^2/2}
\label{eq:pulseG}
\end{equation}
with bandwidth $\pi B = 0.025 \om_o$.

We display the absolute value of the image $\cI(\bx;\bw^\star)$ in the
right plot of Figure \ref{fig:CR}. For comparison, we show in the left
plot of Figure \ref{fig:CR} the image with the uniform weights $w_j =
1/\sqrt{N}$. It has prominent fringes in the cross-range, which are
mitigated by the optimization over the weights. We do not get the best
cross-range resolution with the weights (\ref{eq:IMh8}). The optimal
$\bw^{\rm cr} \in \mathbb{W}$ for focusing in cross-range
has components 
\[
w_j^{\rm cr} = \frac{\beta_j^2(\om_o)}{\|\boldmath{\beta}^2\|_{x_o}},
\qquad j \in \mathbb{J}_{x_o},\qquad \|\boldmath{\beta}^2\|_{x_o} =
\sqrt{\sum_{j\in \mathbb{J}_{x_o}} \beta_j^4(\om_o)}.
\]
It maximizes the ratio of the peak of the image and its mean
square along the cross-range line at $z = 0$,
\begin{equation}
  \cM^{\rm cr}(\bw) = \frac{\left| \cI(\bx_o;\bw)
    \right|^2}{\|\cI((\cdot,0);\bw)\|^2}, \qquad 
  \|\cI((\cdot,0);\bw)\|^2 = \int_0^D dx \, 
  \left|\cI((x,0);\bw)\right|^2,
\end{equation}
and gives the image
\[
\cI(\bx;\bw^{\rm cr}) \approx \frac{1}{4 \|\boldmath{\beta}^2\|_{x_o}}
\sum_{j\in \mathbb{J}_{x_o}} \phi_j(x) \phi_j(x_o) e^{-i
  \beta_j(\om_o) z} f\left[-\beta'_j(\om_o) B z\right].
\]
We show it in the middle plot of Figure \ref{fig:CR}, and indeed, it
has smaller fringes along the axis $z=0$. However, the range resolution is worse than that
given by the optimal weights. 

It is easy to see that the optimal $\bw^{\rm r} \in \mathbb{W}$ for
focusing in range, the maximizer of
\begin{equation}
  \cM^{\rm r}(\bw) = \frac{\left| \cI(\bx_o;\bw)
    \right|^2}{\|\cI((x_o,\cdot);\bw)\|^2}\, , \qquad 
  \|\cI((x_o,\cdot);\bw)\|^2 = \int_{-\infty}^\infty dz \, 
  \left|\cI((x_o,z);\bw)\right|^2,
\end{equation}
has the components 
\[
w_j^{\rm r} = C^{\rm r} \frac{\beta_j(\om_o)}{\phi_j^2(x_o)} \, ,
\qquad j \in \mathbb{J}_{x_o} \, ,
\]
with constant 
\[
C^{\rm r} = 1/\sqrt{ \displaystyle \sum_{l\in \mathbb{J}_{x_o}}
  \beta_l^2(\om_o)/ \phi_l^4(x_o)}.
\]
When $x_o = D/2$ we have $\phi_j^2(x_o) = 1/D$ for all $j$, and
therefore $\bw^{\rm r} = \bw^\star$.  For all other $x_o$ we have
$\bw^{\rm r} \ne \bw^\star$, and the image is given by
\[
\cI(\bx;\bw^{\rm r}) = \frac{C^{\rm r}}{4 } \sum_{j\in
  \mathbb{J}_{x_o}} \frac{ e^{-i \beta_j(\om_o) z}}{\beta_j(\om_o)}
f\left[-\beta'_j(\om_o) B z\right]
\]

Our optimization finds a compromise between cross-range and range
focusing, which is achieved at the maximum of the figure of merit
$\cM(\bw)$.  We can determine explicitly the cross-range and range
resolution of $\cI(\bx;\bw^\star)$ under the assumption that $N \gg
1$ (that is, $D \gg \la_o$).  Then, we can replace the sum over the modes by an integral over
the variable $u = j/N \in (0,1]$, and obtain from the expressions of
$\phi_j$ and $\beta_j$ that
\begin{align}
  \cI((x,0);\bw^\star) &\sim \int_0^1 du\, \frac{\cos \left[ u k_o
      (x-x_o)\right]}{\sqrt{1-u^2}} + \int_0^1 du \,\frac{\cos \left[ u
      k_o (x+x_o)\right]}{\sqrt{1-u^2}} \nonumber \\
  &= \frac{\pi}{2} J_0\left[k_o (x-x_o)\right] + \frac{\pi}{2}
  J_0\left[k_o (x+x_o)\right] \nonumber \\
  & \approx \frac{\pi}{2} J_0\left[k_o (x-x_o)\right],
\end{align}
where $J_0$ is the Bessel function of the first
kind of order zero, and $\sim$ denotes approximate, up to a
multiplicative constant. The cross-range resolution is estimated as
the distance between the peak of $J_0$ that occurs at zero, and its
first zero, that occurs at $k_o|x-x_o| \approx 2.4$. We obtain that
\begin{equation}
|x-x_o| \lesssim \frac{2.4 \la_o}{2 \pi} \sim \frac{\la_o}{2},
\end{equation}
which is basically the diffraction limit of half a wavelength.

\begin{figure}[t]
\begin{center}
\includegraphics[width=4.2cm]{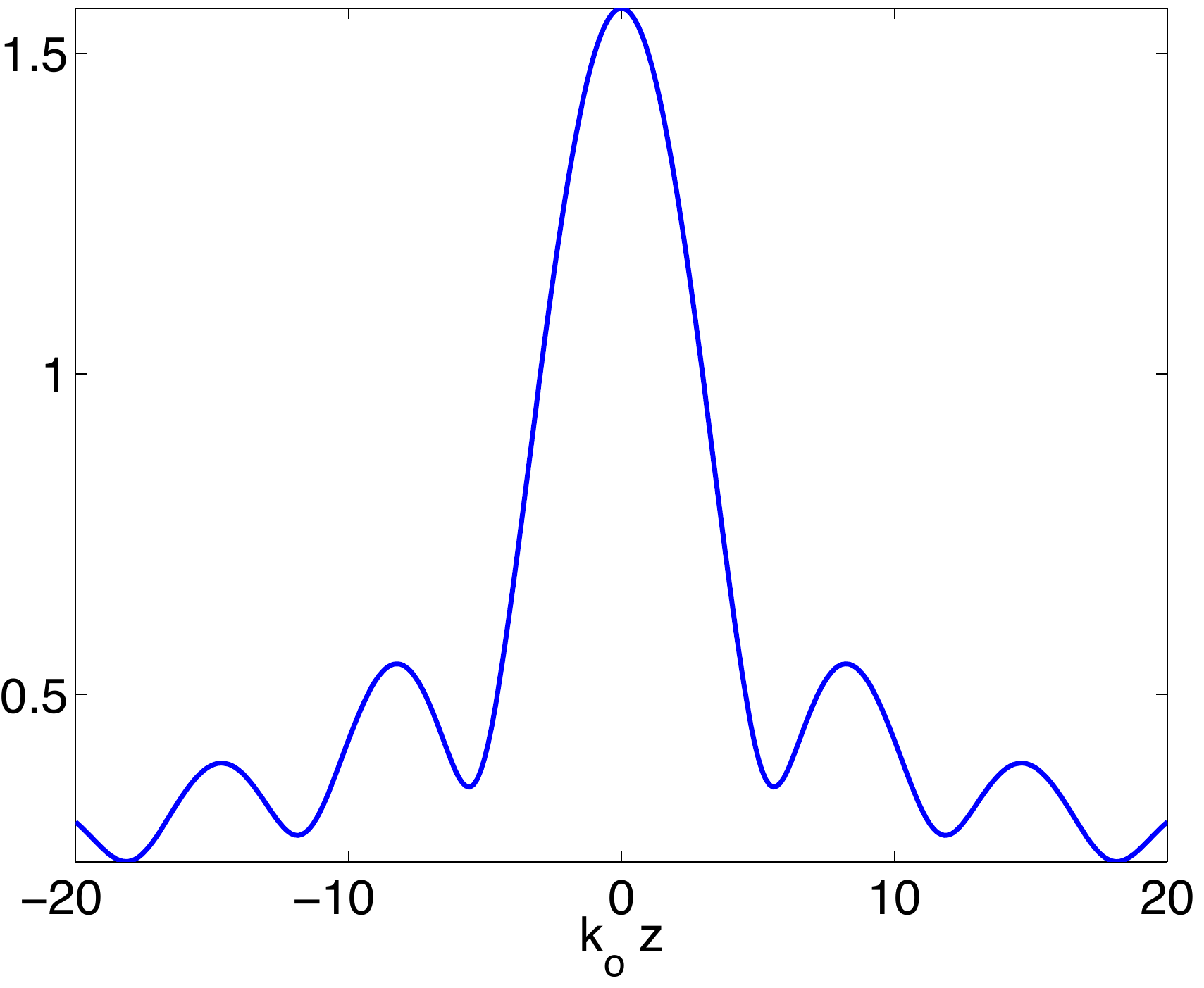}
\end{center}
\caption{The absolute value of the right hand side in
  (\ref{eq:rangeth}), which describes the range resolution, as a
  function of $k_o z$.}
\label{fig:RANGE}
\end{figure}
For the focusing in range we have
\begin{align}
  \cI((x_o,z);\bw^\star) &\sim \int_0^1 du\, \frac{
  e^{-i k_o z \sqrt{1-u^2}}}{\sqrt{1-u^2}}
  \exp\left[-\frac{(B z)^2}{2 c_o^2 (1-u^2)}\right] \nonumber \\
  &\approx \int_0^1 du\, \frac{ e^{-i k_o z
    \sqrt{1-u^2}} }{\sqrt{1-u^2}} \nonumber \\
  &= \frac{\pi}{2} J_o(k_o z) - \frac{i\pi}{2} H_0(k_o z),
\label{eq:rangeth}
\end{align} 
where we used that $x_o = D/2$ and neglected the effect of the pulse
because $B \ll \om_o$.  The result is in terms of the Bessel function
$J_o$ and the Struve function $H_0$, and it is plotted in Figure
\ref{fig:RANGE}. The range resolution is of the order $\la_o$.

\subsection{Random waveguides}
\label{sect:random}
We use the figure of merit (\ref{eq:ThM}) to analyze the optimal
weights for imaging in random waveguides. This is justified as long as
the imaging process remains statistically stable, as explained in
section \ref{sect:WOPT}.  When the data become incoherent, that is
when the array is farther than the scattering mean free path of all
the modes, the weights predicted by the analysis are not useful. The
images have spurious peaks that change unpredictably with the
realization of the random waveguides (the random fluctuations of $\cI$
dominate the mean $\EE[\cI]$). We cannot use coherent imaging for such
data no matter how we weight its components.

\subsubsection{The first two moments of the imaging function}
The imaging function follows from equations (\ref{eq:randp}), 
(\ref{eq:IM3}) and (\ref{eq:ModeDec}), for the full aperture array
\begin{equation}
  \cI(\bx;\bw) \approx \sum_{j=1}^N \int_{-\infty}^\infty \frac{d \om}{2 \pi} \sum_{j=1}^N \frac{w_j \, \overline{\hat a_j(\om,z_\cA)}}{2 
    i \beta_j^{3/2}(\om)} \phi_j(x) e^{-i \beta_j(\om)z}\, . 
\end{equation}
We compute its mean and intensity using the moment formulas
(\ref{eq:meanA}) and (\ref{eq:2ndM}). 
We have
\begin{equation}
\label{eq:meanI}
\EE\left[\cI(\bx;\bw)\right] \approx \frac{1}{4}\sum_{j=1}^N 
\frac{w_j}{\beta_j^2(\om_o)} \phi_j(x) \phi_j(x_o) 
\cF_j(z) \exp\left[-\frac{z_\cA}{\cS_j(\om_o)} - 
i \frac{z_\cA}{\cL_j(\om_o)}\right] \,,
\end{equation}
with mode pulses $\cF_j(z)$ defined in (\ref{eq:IMh3}). This expression
is similar to that of the imaging function in ideal waveguides given
by (\ref{eq:IMh2}), except that the contribution of the $j-$th mode is
damped on the range scale $\cS_j$ and is modulated by oscillation on
the range scale $\cL_j$. This oscillation must be removed in order to
focus the image, which is why we should allow the weights $w_j$ to be
complex.

The intensity of the image is 
\begin{align}
  \EE\left[\left|\cI(\bx ; \bw ) \right|^2\right] \approx &
  \frac{1}{4} \sum_{j,j'=1}^{N} \frac{w_j
    \overline{w_{j'}}}{\beta_j^{3/2}(\om_o)\beta_{j'}^{3/2}(\om_o)}
  \int_{-\infty}^\infty \frac{d \om}{2 \pi} \int_{-\infty}^\infty
  \frac{d \om'}{2 \pi} \EE\left[ \overline{\hat{a}_j(\omega,z_\cA)}\hat
    a_{j'}(\om',z_\cA) \right]
  \nonumber \\
  &\times \, \phi_j(x) \phi_{j'}(x)e^{i\left[\beta_{j'}(\om')-
      \beta_j(\omega)\right] z}
\end{align}
and its square $L^2$ norm is given by 
\begin{align*}
  &\EE\left[ \left\| \cI(\cdot; \bw) \right\|^2\right] = \int_0^D dx
  \int_{-\infty}^\infty dz \EE\left[\left|\cI(\bx ; \bw )
    \right|^2\right] \nonumber \\
  &\hspace{0.3in}=\frac{1}{4} \sum_{j=1}^{N} \frac{|w_j|^2 }{\beta_j^{3}(\om_o)}
  \hspace{-0.02in}\int_{-\infty}^\infty \frac{d \om}{2 \pi}
  \hspace{-0.02in}\int_{-\infty}^\infty \frac{d \om'}{2 \pi} \EE\left[
    \overline{\hat{a}_j(\omega,z_\cA)}\hat a_{j}(\om',z_\cA) \right]
  \hspace{-0.02in} \int_{-\infty}^\infty \hspace{-0.1in}dz \,
  e^{i\left[\beta_{j}(\om')- \beta_j(\omega)\right] z} \\
  &\hspace{0.3in}\approx \frac{1}{4} \sum_{j=1}^{N} \frac{|w_j|^2
  }{\beta_j^{3}(\om_o)\beta'_j(\om_o)} \int_{-\infty}^\infty \frac{d
    \om}{2 \pi}\, \EE\left[ \left|\hat{a}_j(\omega,z_\cA)\right|^2
  \right]\, ,
\end{align*}
because of the orthonormality of the eigenfunctions $\phi_j(x)$.
Recalling the moment formula (\ref{eq:2ndM}) and using equation 
(\ref{eq:ModeSpeed}), we obtain 
\begin{align}
  \EE\left[ \left\| \cI(\cdot; \bw) \right\|^2\right] &\approx
  \frac{c_o \|f\|^2}{16 k_o B} \sum_{j=1}^N
  \frac{|w_j|^2}{\beta_j^2(\om_o)} \sum_{l=1}^N
  \frac{\phi_l^2(x_o) T_{jl}(\om_o,z_\cA)}{\beta_l(\om_o) }.
\end{align}

\subsubsection{Optimal weights}
The weights must compensate for the oscillations in (\ref{eq:meanI}) in
order for $\EE\left[\cI(\bx;\bw)\right]$ to peak at the source
location $\bx_o$. Thus, we let 
\begin{equation}
w_j = w_j^+ \exp\left[i \frac{z_\cA}{\cL_{j}(\om_o)}\right], \qquad 
w_j^+ = |w_j|,
\label{eq:weightsmod}
\end{equation}
and maximize
\begin{align}
  \cM_{\rm th}(\bw^+) &= \frac{ \left| \EE \left[ \cI(\bx_o;\bw)
      \right] \right|^2 }{ \EE \left[ \| \cI(\cdot;\bw) \|^2 \right]}
  \sim \frac{ \left[ \displaystyle \sum_{j=1}^N \frac{w_j^+
        \phi_j^2(x_o)}{\beta_j^2(\om_o)}
      \exp\left(-\frac{z_\cA}{\cS_j(\om_o)}\right)
    \right]^2}{\displaystyle \sum_{j=1}^N
    \frac{(w_j^+)^2}{\beta_j^2(\om_o)}\displaystyle \sum_{l=1}^N
    \frac{\phi_l^2(x_o) \, T_{jl}(\om_o,z_\cA)}{\beta_l(\om_o) }},
\label{eq:ThMT}
\end{align}
over the vectors $\bw^+ = (w_1^+, \ldots, w_N^+)^T$ with non-negative
entries, and Euclidian norm $\| \bw^+ \| = 1$. The symbol $\sim$ denotes
approximate, up to a multiplicative constant, as before.

The optimal weights are given by 
\begin{equation}
  w_j^+ = \frac{C \, \phi_j^2(x_o) \exp \left(-\frac{z_\cA}{\cS_j(\om_o)}\right)}{
    \displaystyle \sum_{l=1}^N \frac{\phi_l^2(x_o)
     \, T_{jl}(\om_o,z_\cA)}{\beta_l(\om_o) }}, \quad 
  j \in \mathbb{J}_{x_o}\, ,
\end{equation}
with positive constant $C$ determined by the normalization $\|\bw^+\|
= 1$. They are damped exponentially with range on the scale given by
the mode dependent scattering mean free paths $\cS_j$. The
optimization detects the modes that are incoherent, i.e., the indexes
$j$ for which $z_\cA > \cS_j(\om_o)$, and suppresses them in the data.


\section{Numerical simulations}
\label{sect:num}
In this section we present numerical simulations and compare the
results with those predicted by the theory. The setup is as described
in section \ref{sect:scat}, with autocorrelation functions
(\ref{eq:Rnu}) and (\ref{eq:RmuG}) of the perturbations of the
boundary and of the wave speed, in a waveguide of depth $D = 20\la_o$.
All lengths are scaled by the central wavelength $\la_o$, and the
bandwidth satisfies $\pi B = 0.0625 \om_o$. For example, we could have
the central frequency $1$kHz and the unperturbed wave speed $c_o =
1$km/s, so that $\la_o = 1$m and $B = 0.125$kHz.  
To illustrate the cumulative scattering effect on the imaging process,
we consider several ranges $z_\cA$ of the array, from $25\la_o$ to
$150 \la_o$. The details on the numerical simulations of the array 
data are in appendix \ref{sect:ap}

We begin in Figure \ref{fig:homo} with the results in an ideal
waveguide, with array at range $z_\cA = 100\la_o$. We plot on the left
the image with the optimal weights and on the right the theoretical
weights (\ref{eq:IMh8}) (in red) and the numerically computed weights
(in blue). The weights are computed by minimizing $1/\cM(\bw)$, with
$\cM$ defined in (\ref{eq:FIG_MERRIT}). The optimization is done with
the MATLAB function \emph{fmincon}, over weights $\bw = (w_1, \ldots,
w_N)^T \in \mathbb{R}^N$, with constraints $w_j \ge 0$ for $j = 1,
\ldots, N$, and normalization $\|\bw\|^2 = 1.$ The image is very
similar to that predicted by the theory (the right plot in Figure
\ref{fig:CR}), and the optimal weights are in agreement, as well.
 
\begin{figure}[t]
  \includegraphics[width=6cm]{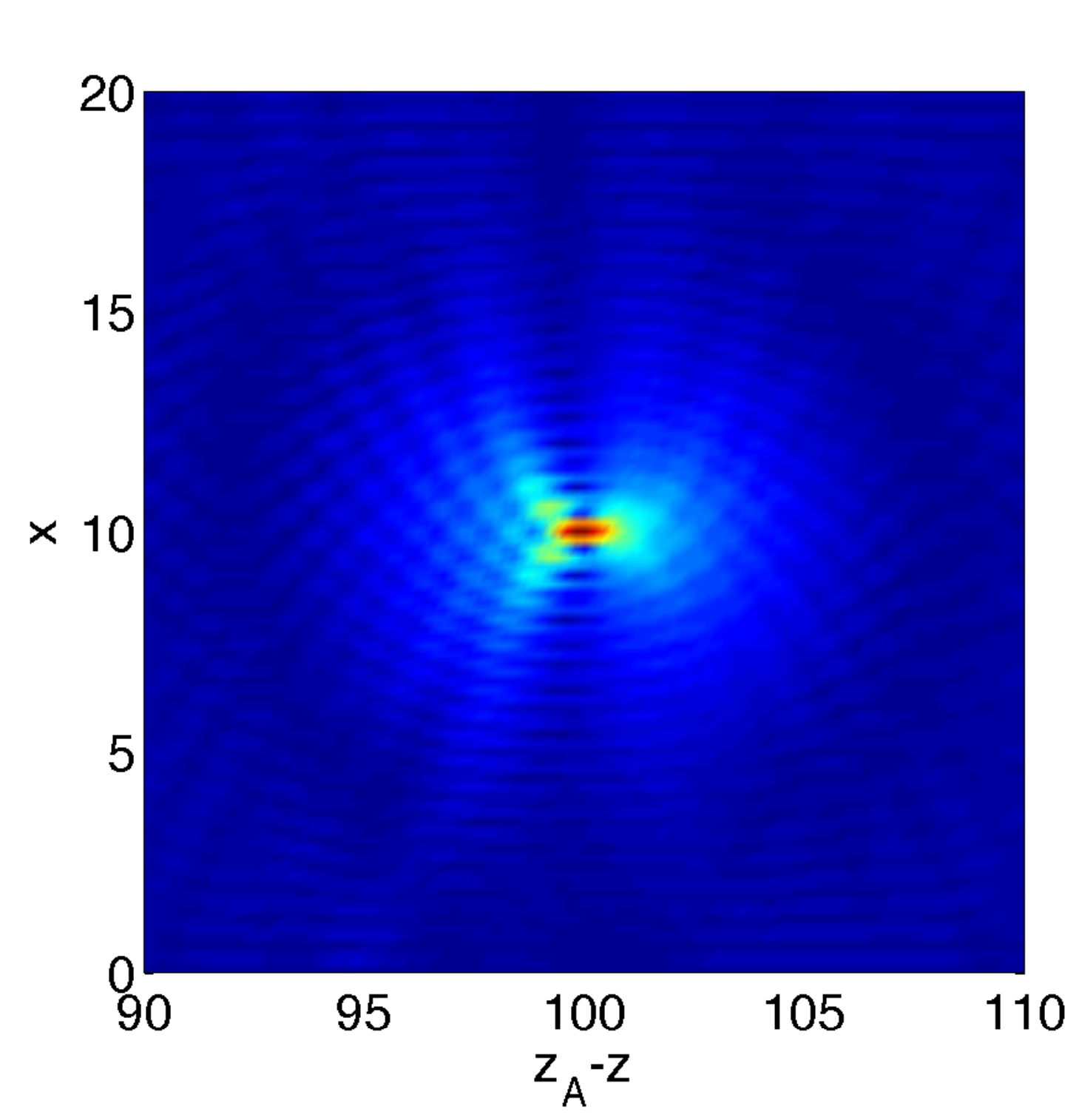}
  \raisebox{0.2in}{\includegraphics[width=6.5cm]{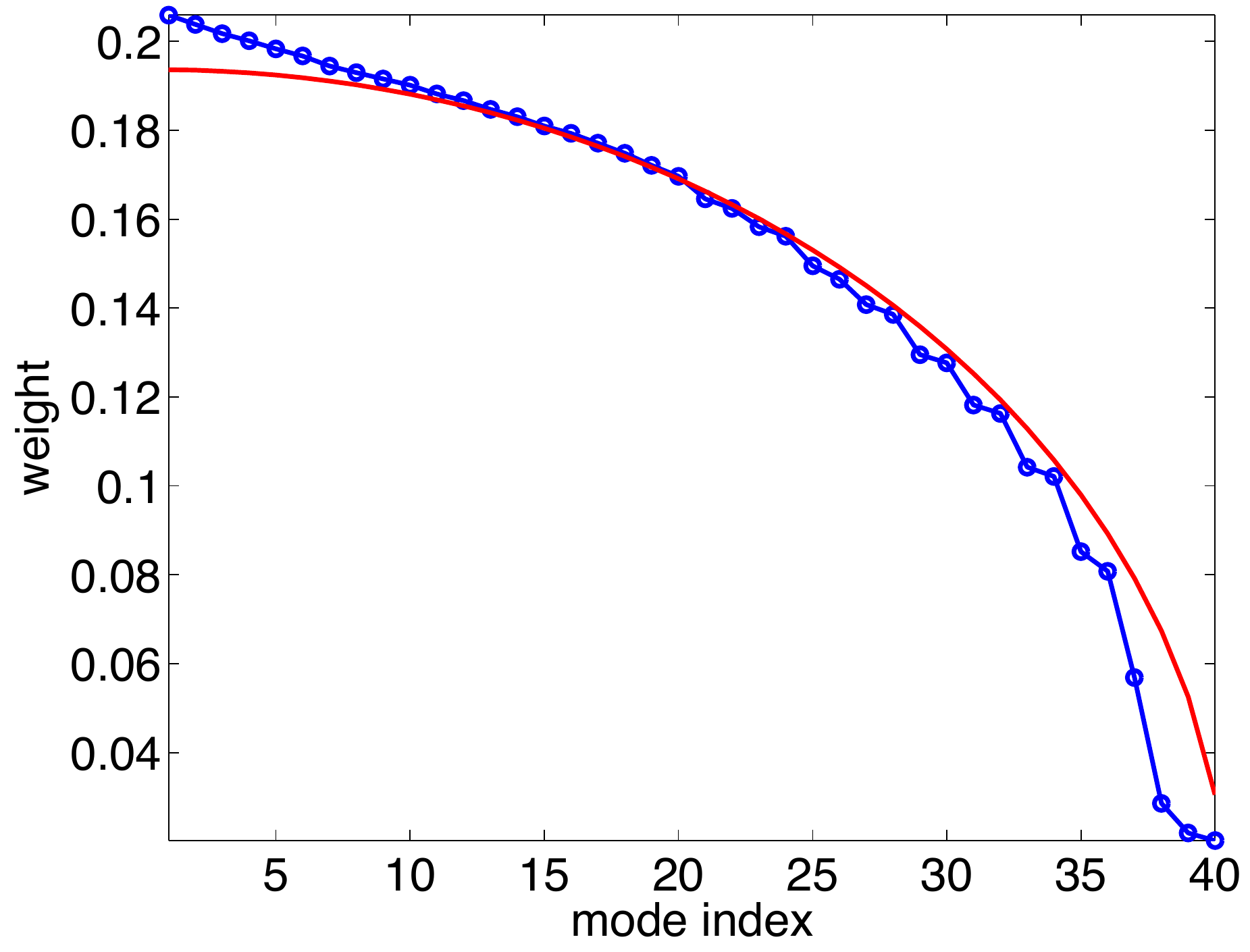}}
  \caption{Homogeneous medium and array range $z_\cA =100 \lambda_o$.
    Left: Image with the numerically computed weights. The abscissa is
    $z_\cA-z$ in $\la_o$ and the ordinate is the cross-range $x$ in
    $\la_o$. Right: Theoretical weights (in red) and numerical ones
    (in blue) vs. mode index.}
\label{fig:homo}
\end{figure}

The analysis for random waveguides in section \ref{sect:random} is
based on the theoretical figure of merit (\ref{eq:ThM}), which is
close to $\cM(\bw)$ only when the image is statistically stable.  The
theory in \cite{ABG-12,kohler77,book07} predicts that stability holds
for the given bandwidth, in the asymptotic limit $\eps \to 0$.  We
have a finite $\eps$, and to stabilize the optimization so that we can
compare it with the theory, we need to work with a slight modification
of the figure of merit (\ref{eq:FIG_MERRIT}),
\begin{equation}
  \cM_{\rm num}(\bw) = \frac{\left| \left< \cI(\bx_o ; \bw) \right>
  \right|^2}{ \|\cI(\cdot;\bw)\|^2},
\label{eq:FIG_NUM}
\end{equation}
where $\left< \cI(\bx_o ; \bw) \right>$ is a local spatial average of
the image around $\bx_o$.

In our regime the theory predicts that $\cL_j > z_\cA$ for all the
modes that remain coherent, as shown in Figures \ref{fig:1} and
\ref{fig:2}.  Therefore, we can neglect the phase factors in
(\ref{eq:weightsmod}), and optimize directly over positive weights.
The optimization is done with the MATLAB function \emph{fmincon}, as
before, but we regularize it by asking that the weights be monotone
decreasing with the mode index. That is to say, we work with the
constraints
\[
w_{j} \ge w_{j+1} \ge 0, \quad j = 1, \ldots, N-1, \quad 
\sum_{j=1}^N w_j^2 = 1\,.
\]
 
\begin{figure}[t]
\includegraphics[width=6cm]{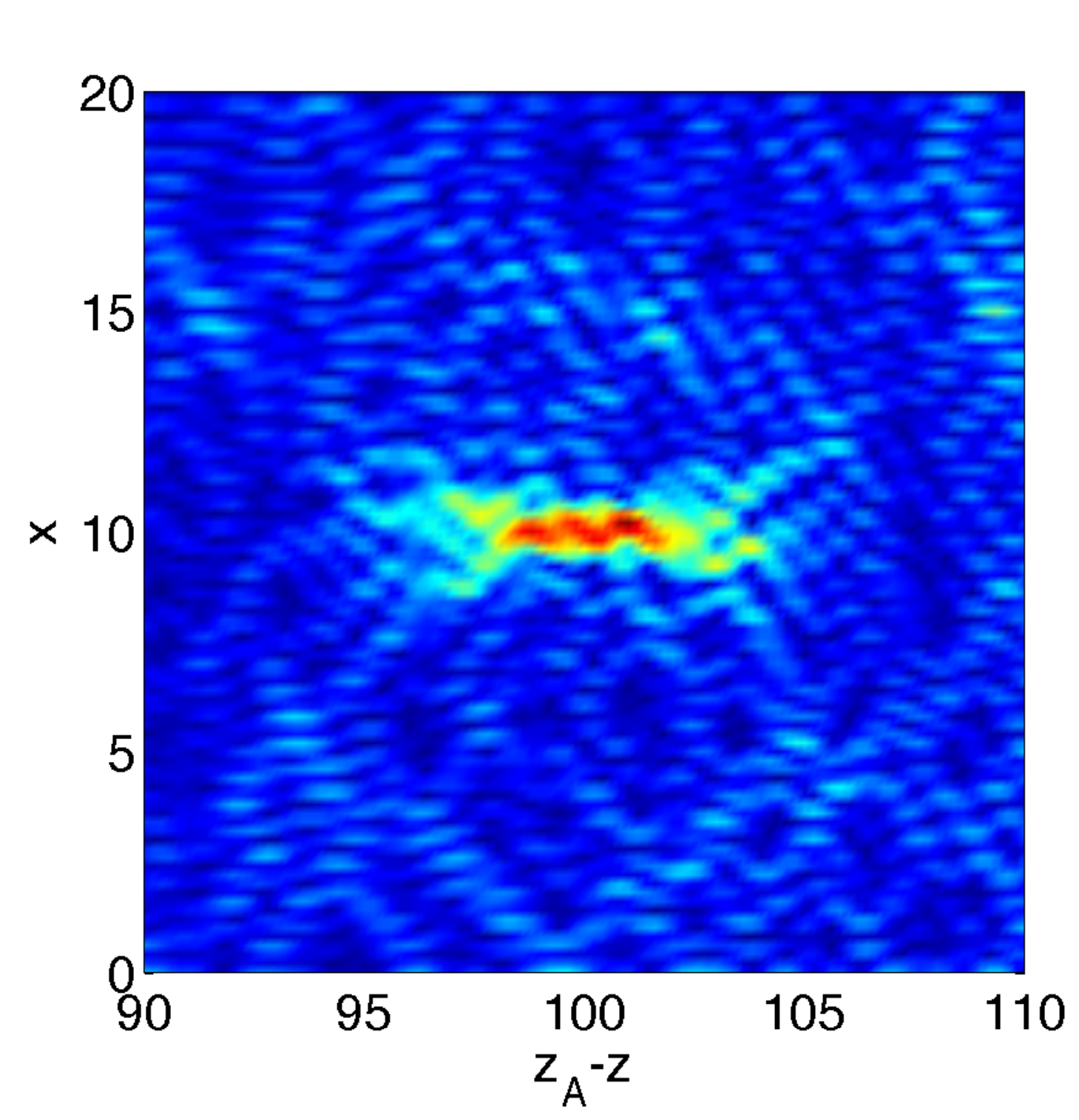}
\includegraphics[width=6cm]{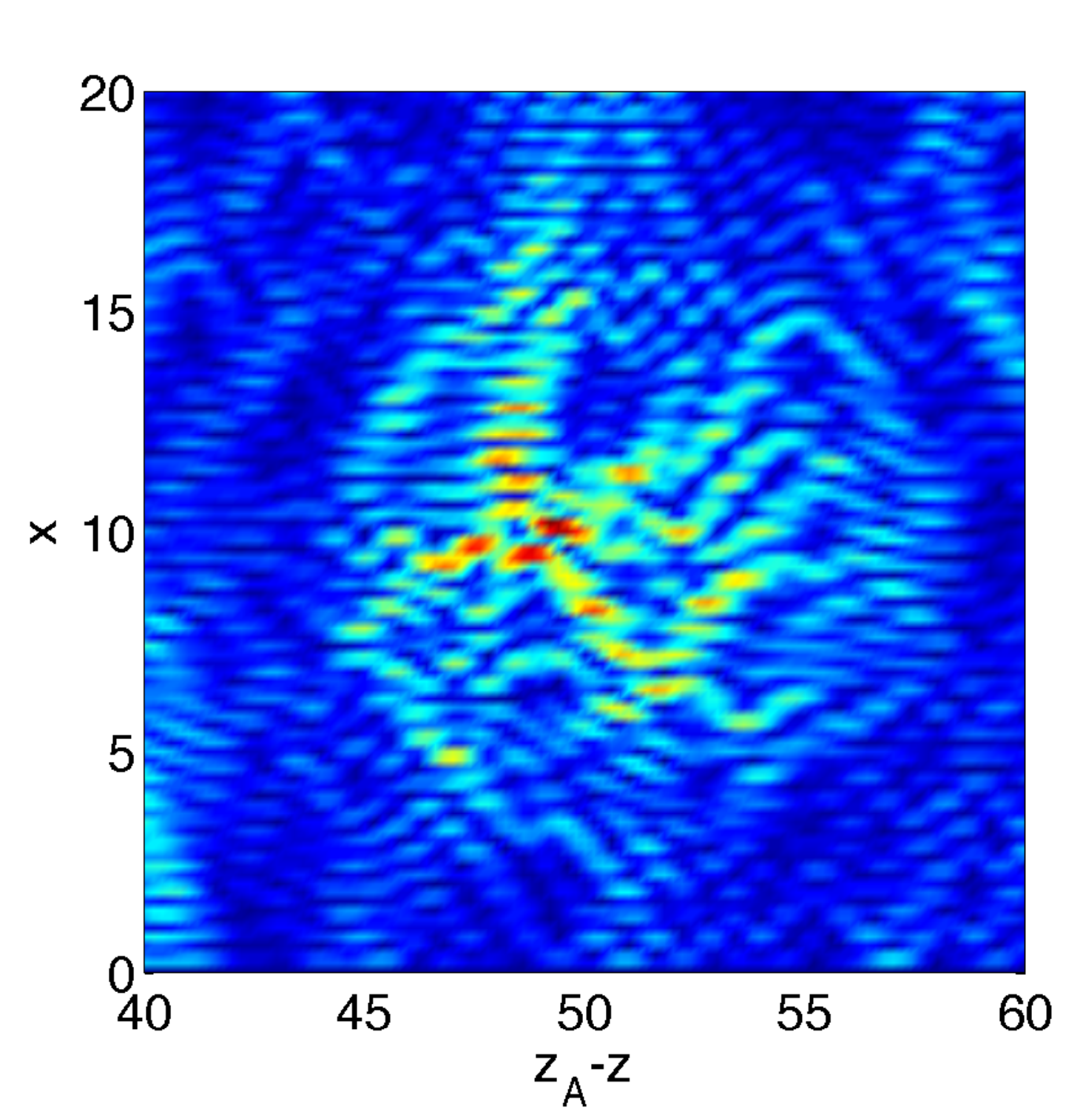}
\caption{Image $\cI(\bx;\bw)$ in waveguide with perturbed boundary and
  array at range $z_\cA = 100 \la_o$ (left) and in waveguide with
  perturbed medium and array at range $z_\cA = 50 \la_o$ (right). The
  weights are uniform $w_j = 1/\sqrt{N}$, for $j = 1, \ldots, N$. The
  abscissa is $z_\cA-z$ in $\la_o$ and the ordinate is the cross-range
  $x$ in $\la_o$. }
\label{fig:KMimages}
\end{figure}

Without weight optimization the images are noisy, with spurious peaks.
We illustrate this in Figure \ref{fig:KMimages}, where we plot
$\cI(\bx;\bw)$ with uniform weights $w_j = 1/\sqrt{N}$, for $j = 1,
\ldots, N$. The image in the left plot is in a waveguide with
perturbed boundary and array at range $z_\cA = 100\la_o$. The image in
the right plot is in a waveguide with perturbed medium and array at
range $z_\cA = 50 \la_o$. Both images are noisy.
The results in Figure \ref{fig:1} predict that
half of the modes remain coherent at $z_\cA = 100\la_o$ in the
waveguide with perturbed boundaries ($\cS_j > 100 \la_o$ for $j = 1,
\ldots, N/2$). Therefore the image is not bad, and can be improved
further by the optimization, as shown below.  The results in Figure
\ref{fig:2} show that all the modes are almost incoherent at $z_\cA =
50\la_o$ in the waveguide with perturbed medium ($\cS_j < 70 \la_o$
for $j = 1, \ldots, N$). The image is noisy, with prominent
spurious peaks, and cannot be improved by optimization, as shown
below.

\begin{figure}[t]
  \includegraphics[width=6cm]{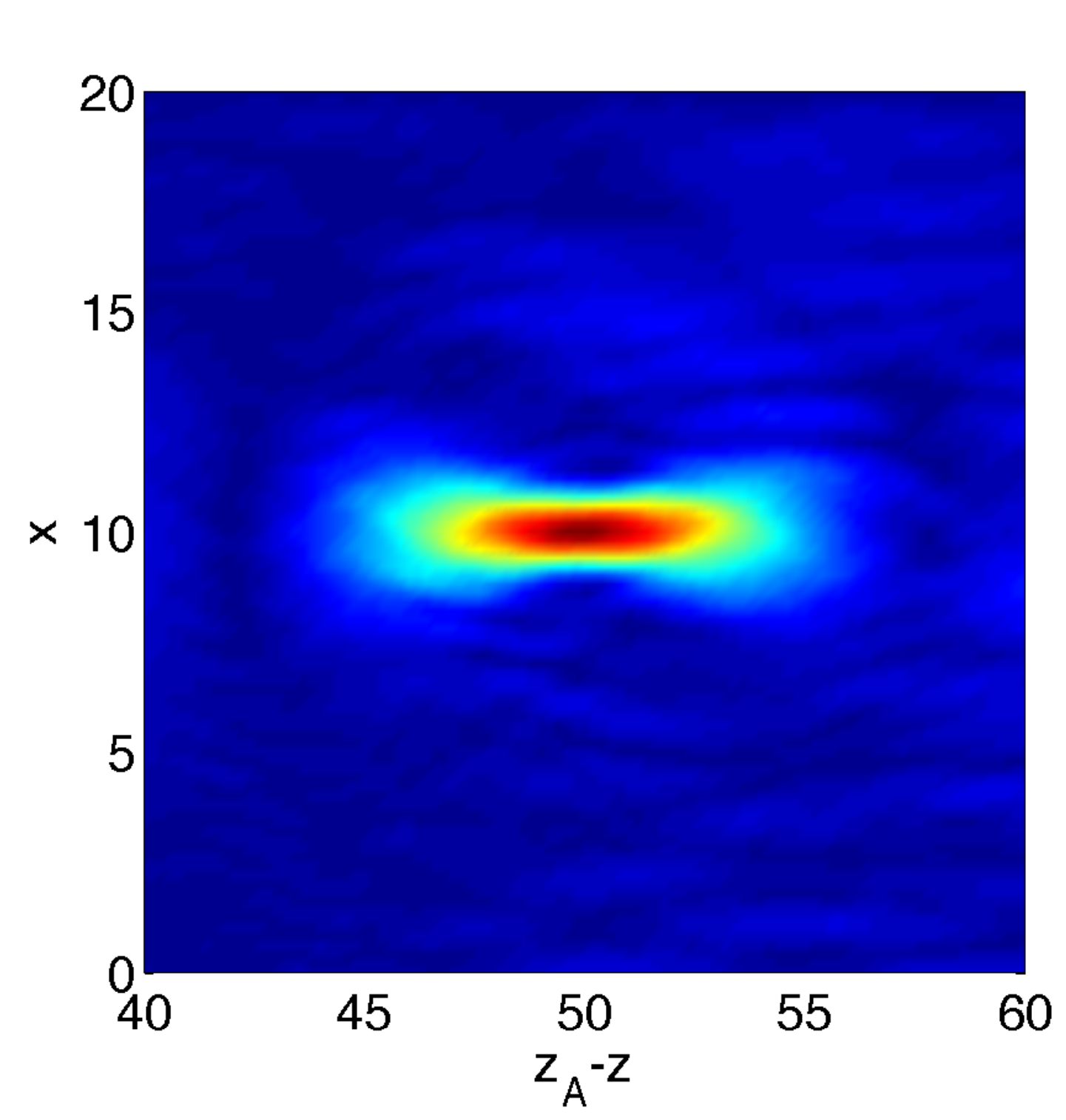}
  \raisebox{0.2in}{\includegraphics[width=6.5cm]{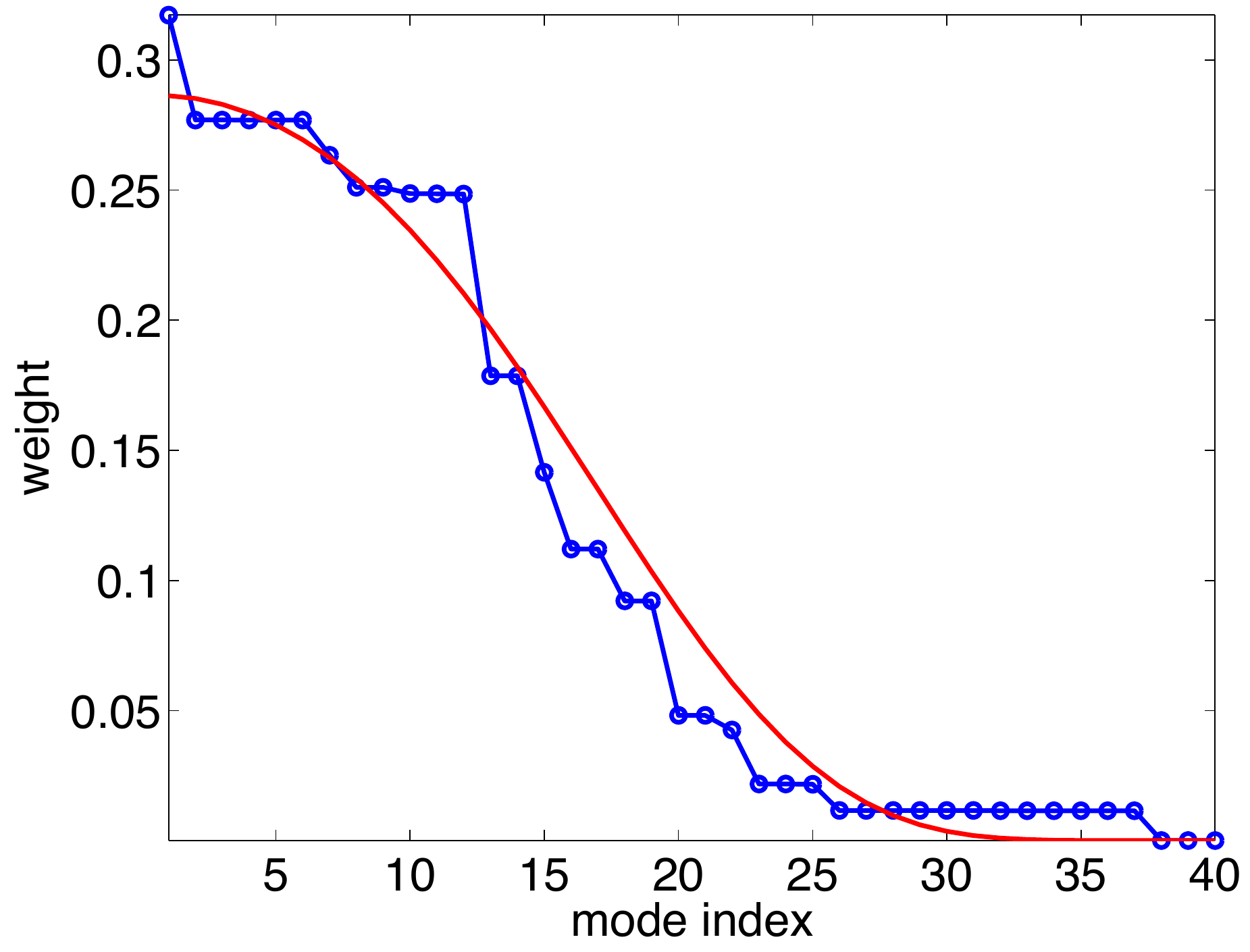}}
  \caption{Waveguide with perturbed boundary and array range $z_\cA
    =50\lambda_o$.  Left: Image with the numerically computed
    weights. The abscissa is $z_\cA-z$ in $\la_o$ and the ordinate is
    the cross-range $x$ in $\la_o$. Right: Theoretical weights (in
    red) and numerical ones (in blue) vs. mode index.}
\label{fig:50BDRY}
\end{figure}

\begin{figure}[h]
  \includegraphics[width=6cm]{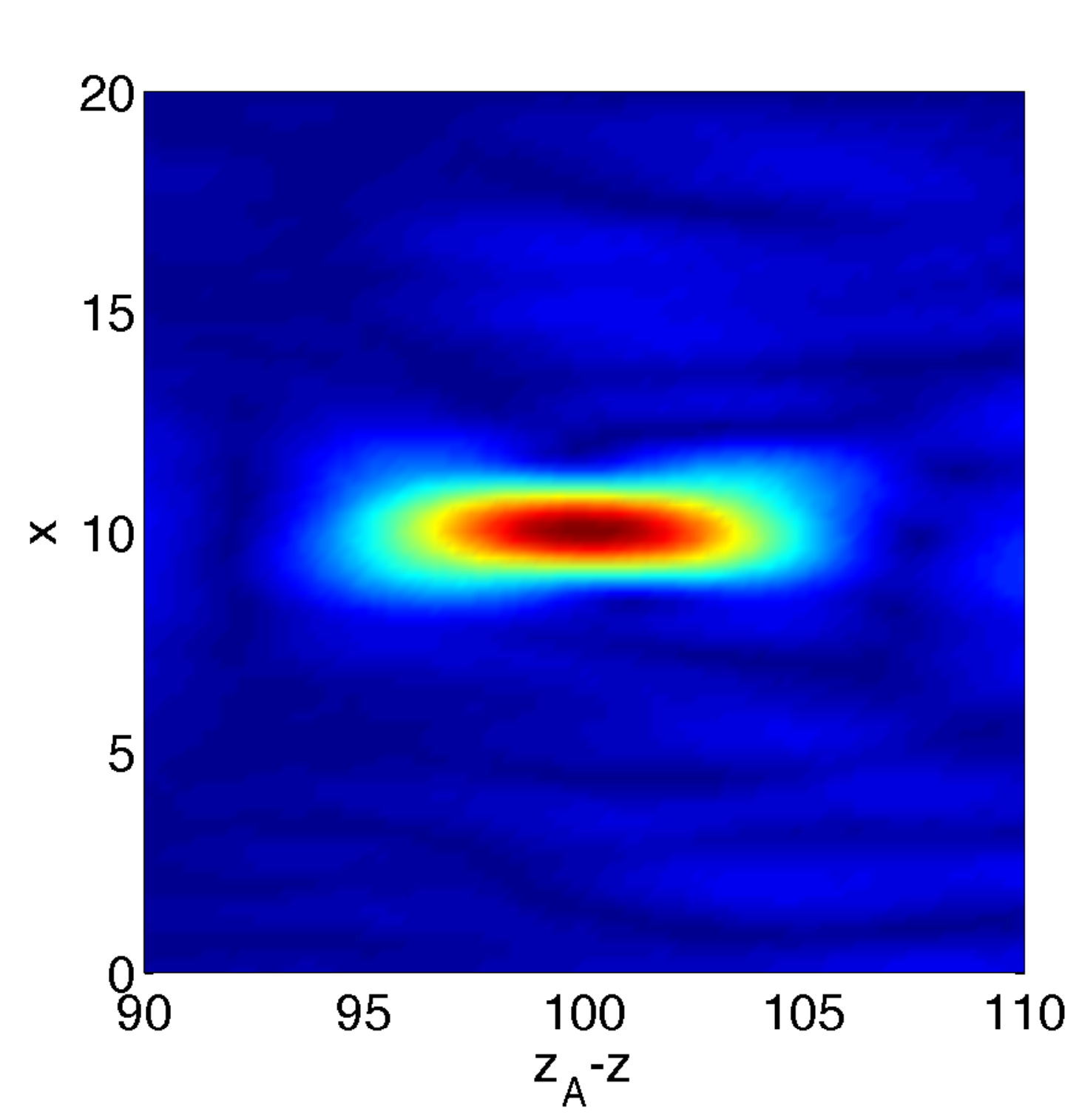}
  \raisebox{0.2in}{\includegraphics[width=6.5cm]{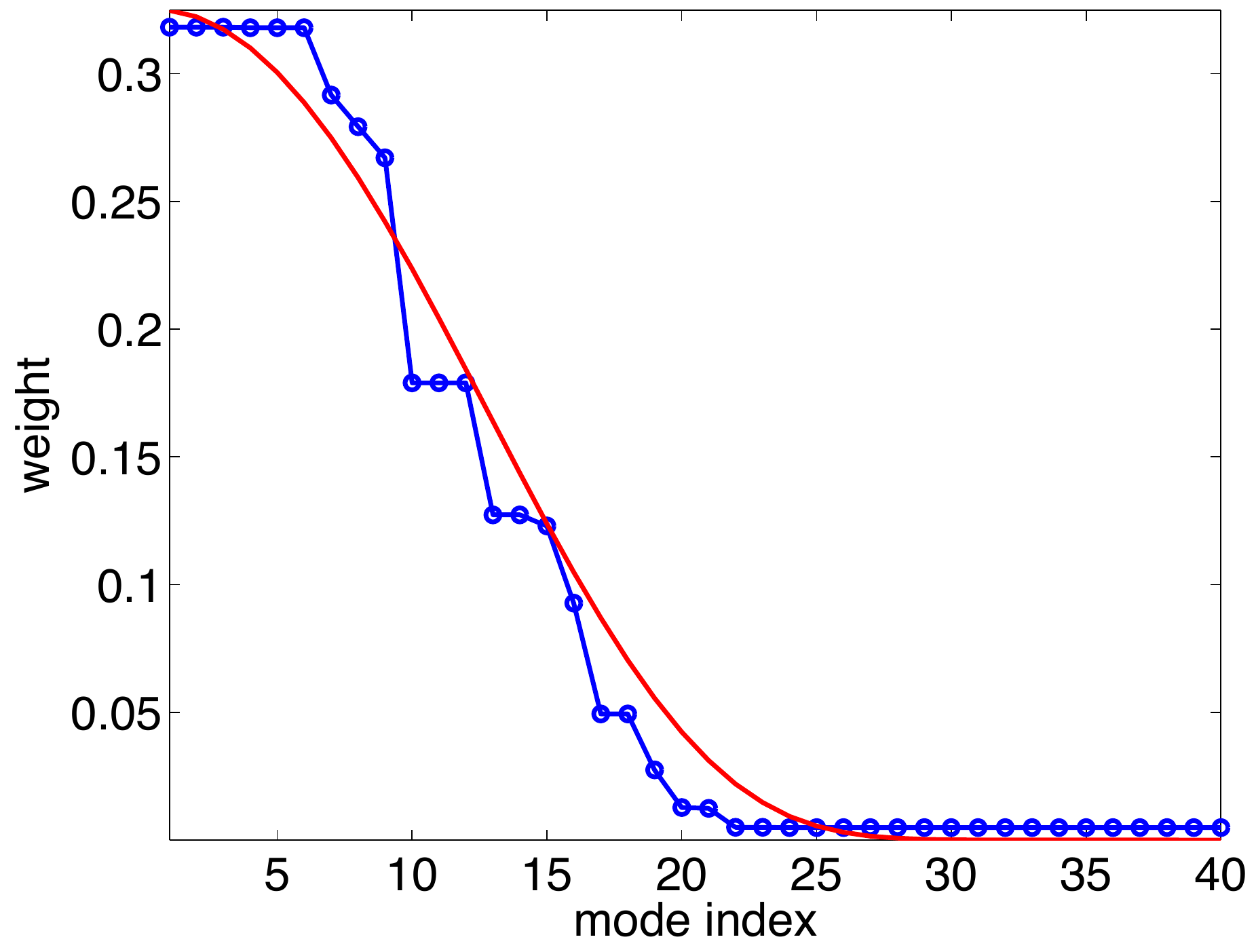}}
  \caption{Waveguide with perturbed boundary and array range $z_\cA
    =100\lambda_o$.  Left: Image with the numerically computed
    weights. The abscissa is $z_\cA-z$ in $\la_o$ and the ordinate is
    the cross-range $x$ in $\la_o$. Right: Theoretical weights (in
    red) and numerical ones (in blue) vs. mode index.}
\label{fig:100BDRY}
\end{figure}

\begin{figure}[h]
  \includegraphics[width=6cm]{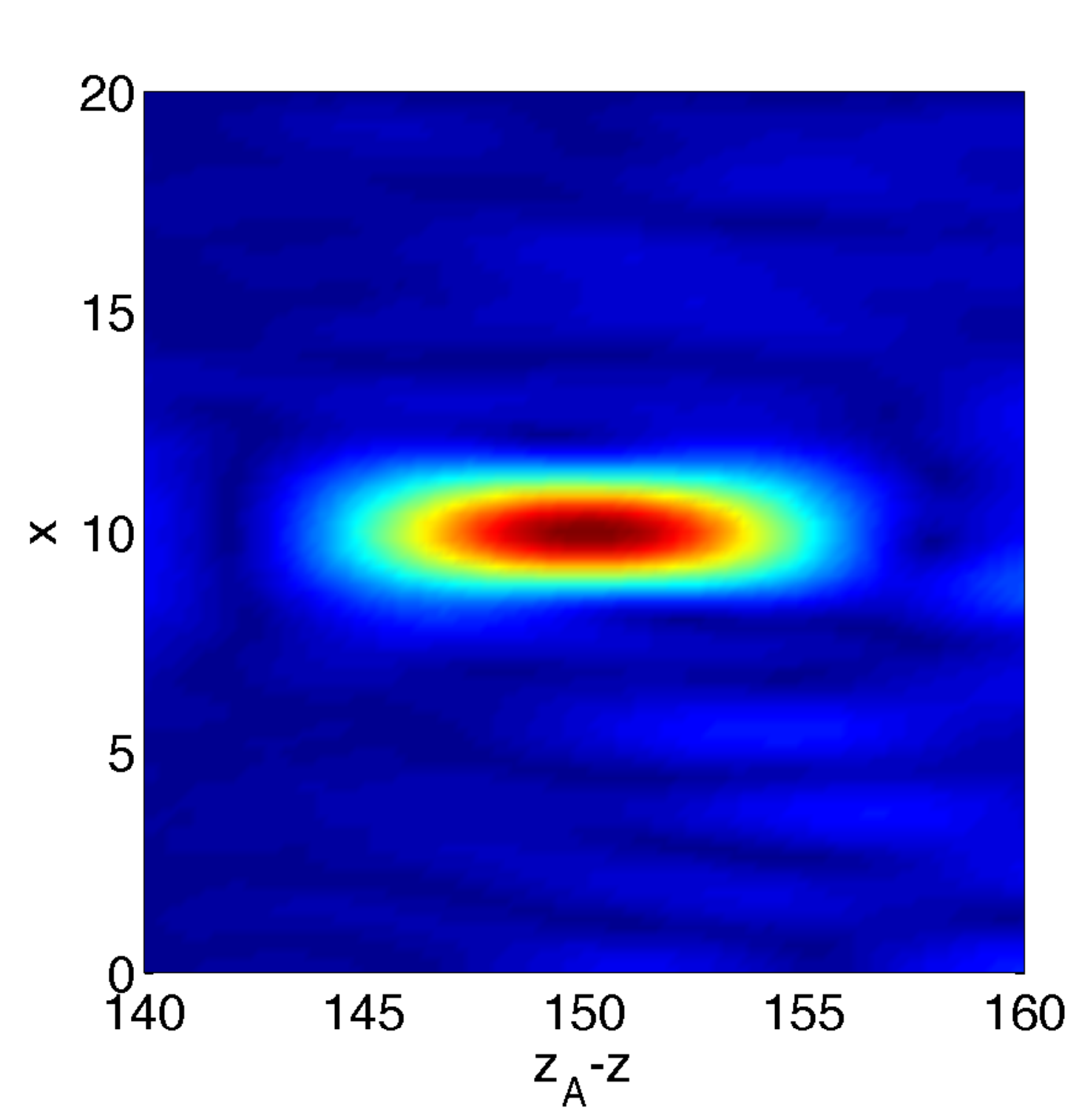}
  \raisebox{0.2in}{\includegraphics[width=6.5cm]{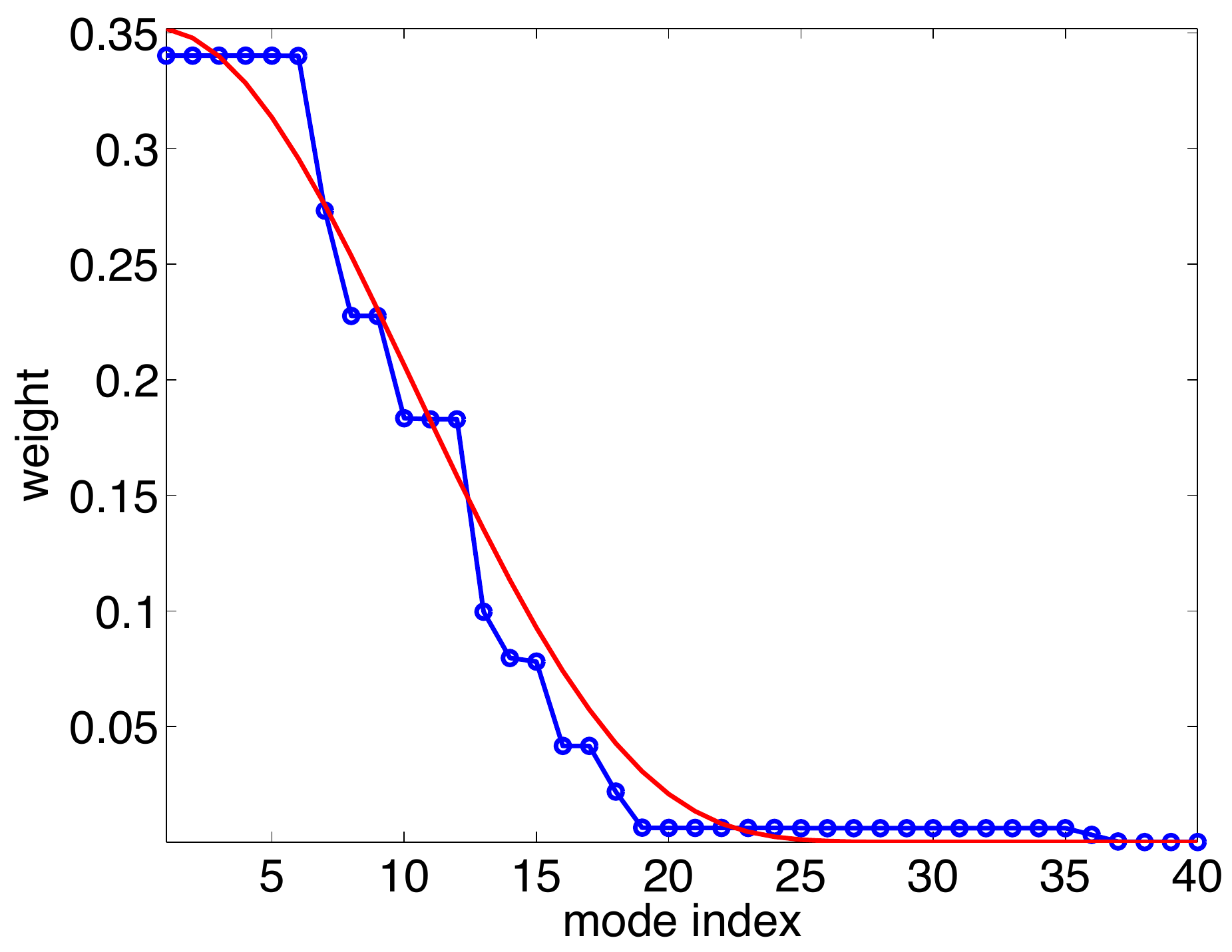}}
  \caption{Waveguide with perturbed boundary and array range $z_\cA
    =150\lambda_o$.  Left: Image with the numerically computed
    weights. The abscissa is $z_\cA-z$ in $\la_o$ and the ordinate is
    the cross-range $x$ in $\la_o$. Right: Theoretical weights (in
    red) and numerical ones (in blue) vs. mode index.}
\label{fig:150BDRY}
\end{figure}

\begin{figure}[h]
  \includegraphics[width=6cm]{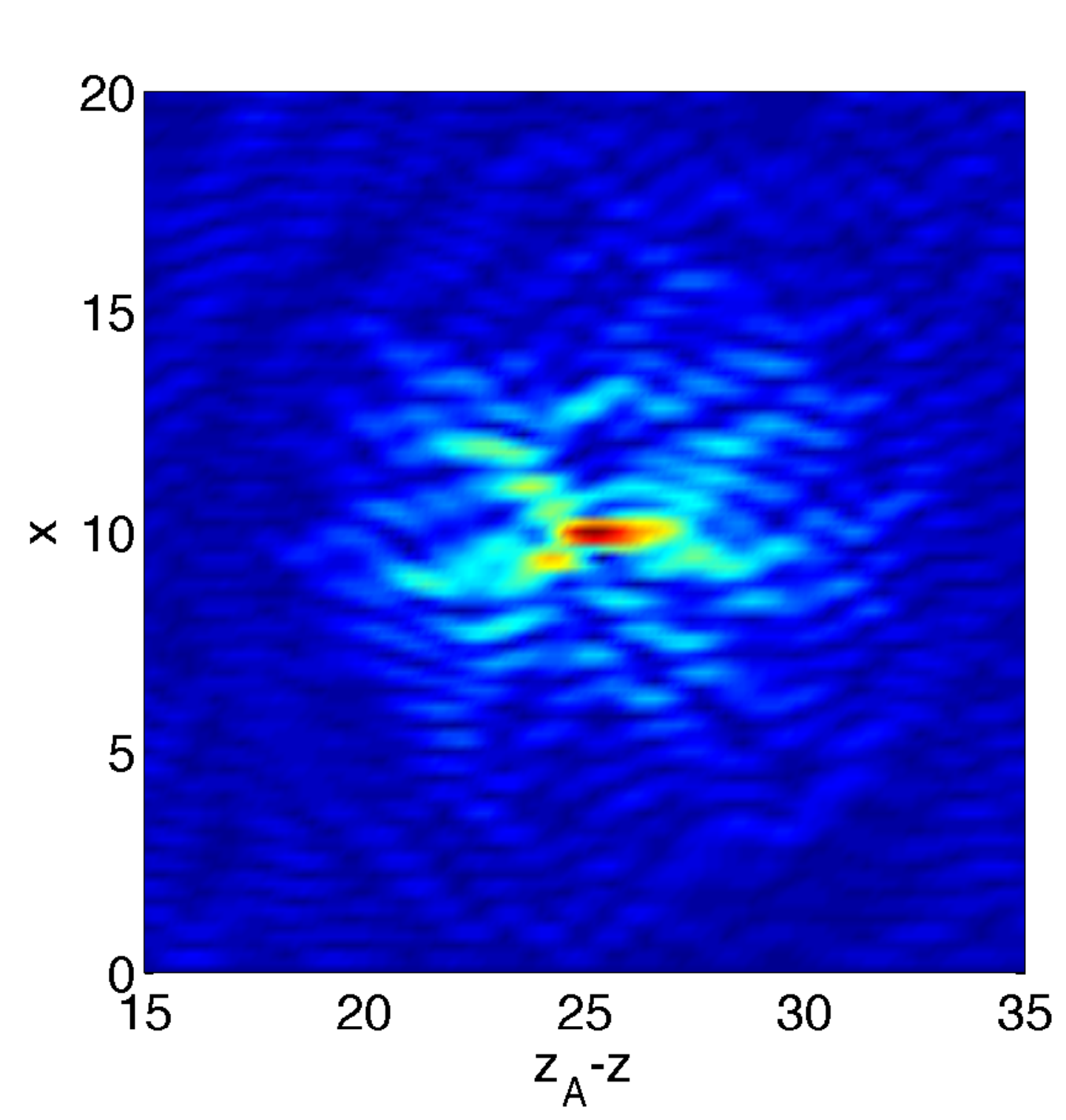}
  \raisebox{0.2in}{\includegraphics[width=6.5cm]{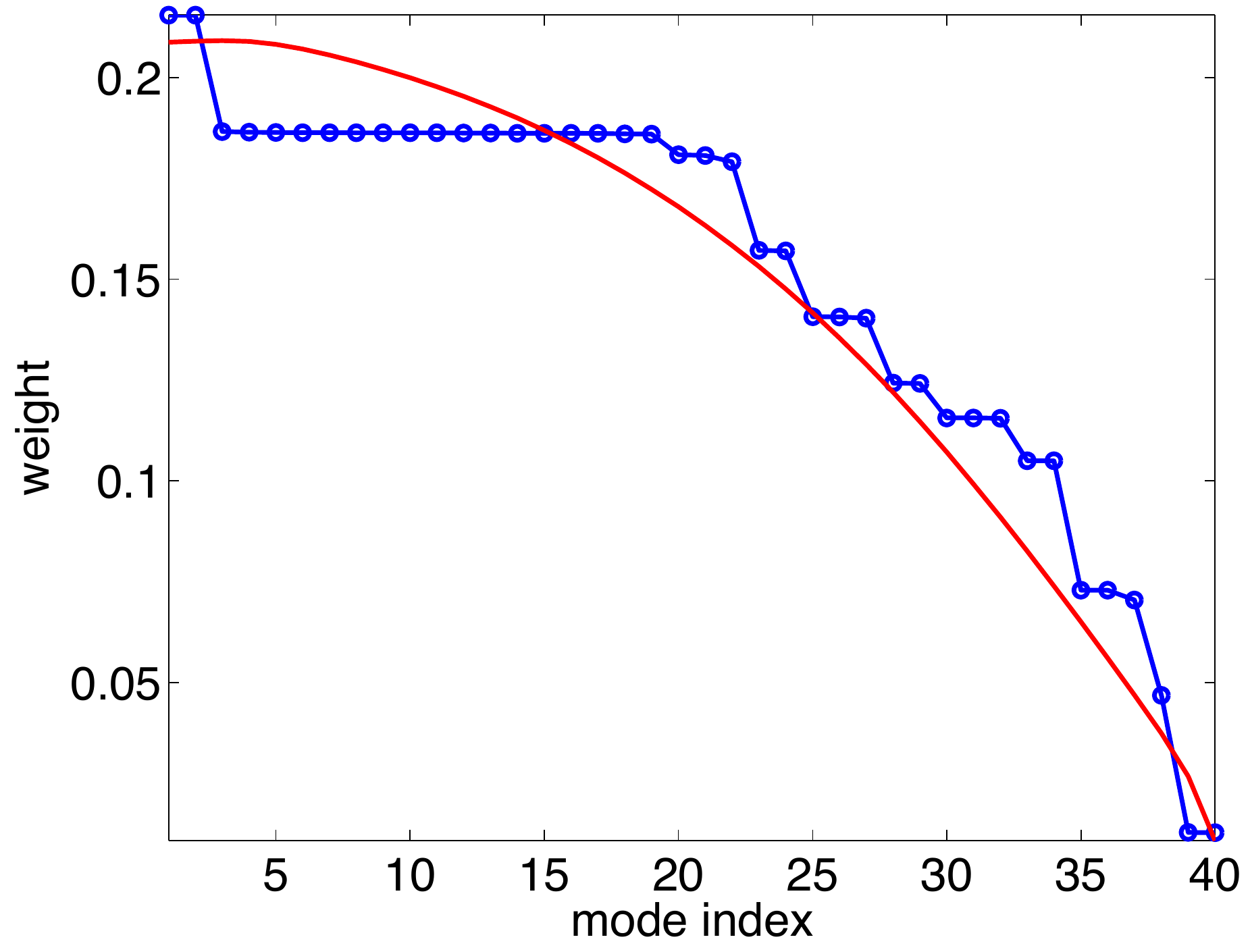}}
  \caption{Waveguide with perturbed medium and array range $z_\cA
    =25\lambda_o$.  Left: Image with the numerically computed
    weights. The abscissa is $z_\cA-z$ in $\la_o$ and the ordinate is
    the cross-range $x$ in $\la_o$. Right: Theoretical weights (in
    red) and numerical ones (in blue) vs. mode index.}
\label{fig:25INT}
\end{figure}

\begin{figure}[h]
  \includegraphics[width=6cm]{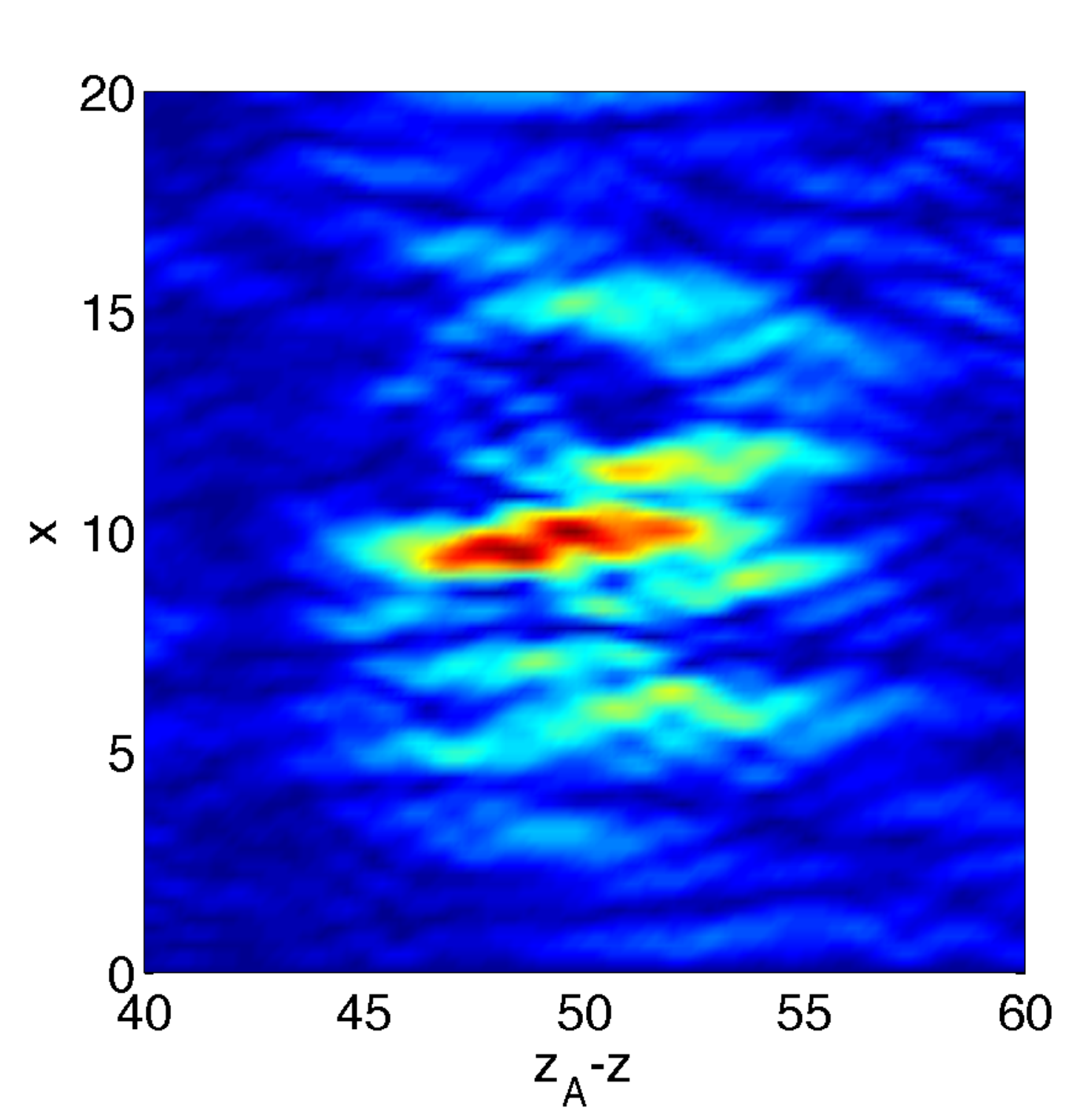}
  \raisebox{0.2in}{\includegraphics[width=6.5cm]{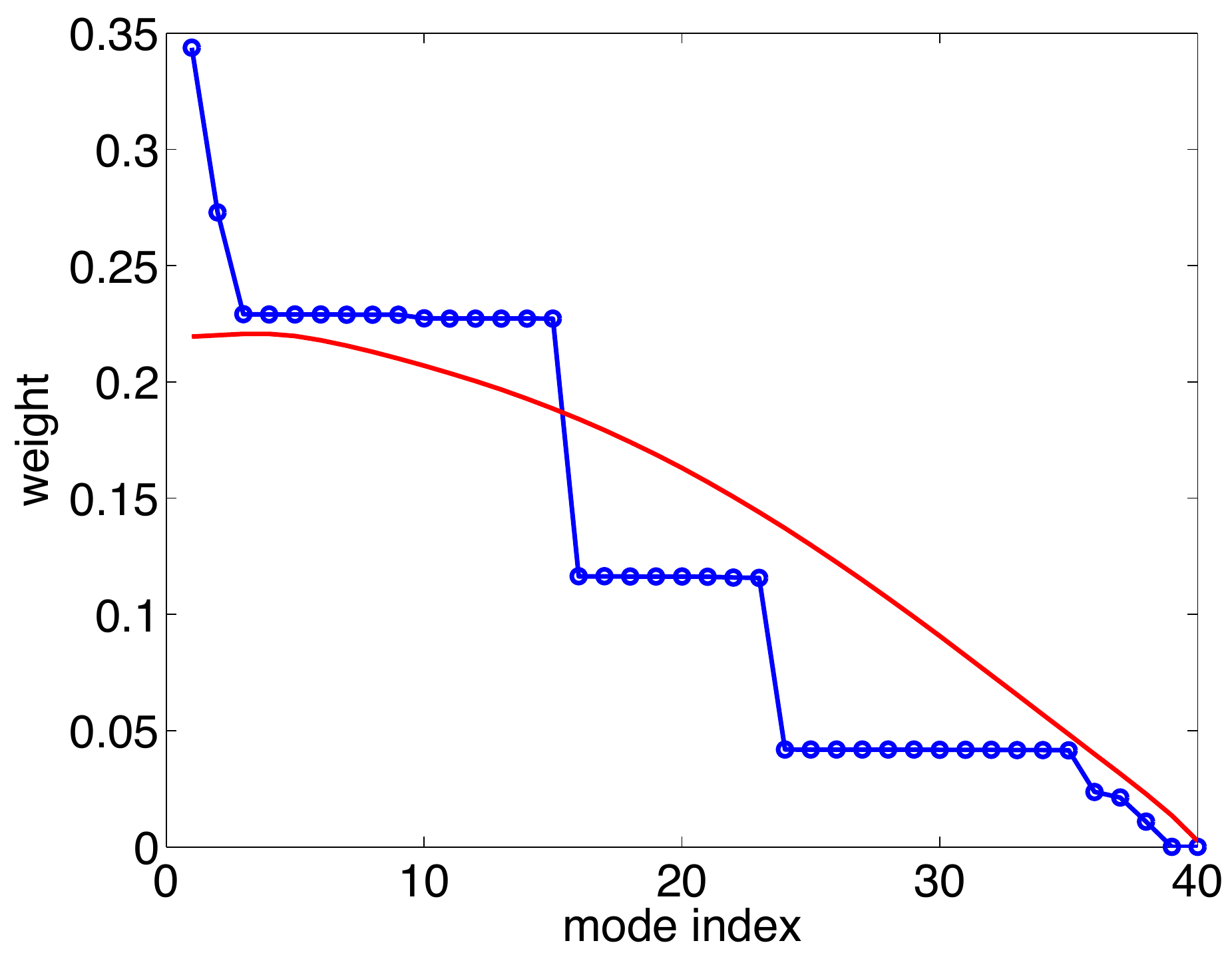}}
  \caption{Waveguide with perturbed medium and array range $z_\cA
    =50\lambda_o$.  Left: Image with the numerically computed
    weights. The abscissa is $z_\cA-z$ in $\la_o$ and the ordinate is
    the cross-range $x$ in $\la_o$. Right: Theoretical weights (in
    red) and numerical ones (in blue) vs. mode index.}
\label{fig:50INT}
\end{figure}

We show in Figures \ref{fig:50BDRY}-\ref{fig:150BDRY} the results of
the optimization in a waveguide with perturbed boundary and array at
ranges $z_\cA = 50 \la_o$, $100 \la_o$ and $150 \la_o$. The local
average of the image in (\ref{eq:FIG_NUM}) is over an interval of
length $\la_o$ in range and of length $\la_o$, $1.5\la_o$ and $2
\la_o$ in cross-range, respectively. The weights obtained with the
numerical optimization are in reasonable agreement with those
predicted by the theory. The resolution of the images deteriorates as
we increase $z_\cA$ because more of the higher indexed modes become
incoherent. 

Figures \ref{fig:25INT}-\ref{fig:50INT} show the results in a
waveguide with perturbed medium and array at ranges $z_\cA = 25\la_o$
and $50 \la_o$. Here there is no trade-off between resolution and
robustness of the image, because most modes lose coherence on roughly
the same range scale. Coherent imaging can be done at range $z_\cA =
25\la_o$, and the numerical weights agree with those predicted by the
theory. However, at range $z_\cA = 50 \la_o$ the optimization fails to 
improve the image.

\section{Summary}
\label{sect:sum}
We have carried out a comparative theoretical and numerical study of
wave scattering in two types of random waveguides with bounded
cross-section: waveguides with random inhomogeneities in the bulk
medium and waveguides with random perturbations of the boundary. The
wave field is a superposition of waveguide modes with random
amplitudes. Coherent imaging relies on the coherent part of the
amplitudes, their expectation. However, this decays with the distance
of propagation due to cumulative scattering at the random
inhomogeneities and boundary perturbations. The incoherent part of the
amplitudes, the random fluctuations gain strength and become dominant
at long ranges.

The characteristic range scales of decay of the coherent part of the
mode amplitudes are called scattering mean free paths. They are
frequency and mode-dependent, and they decrease monotonically with the
mode index.  In waveguides with random boundaries the mode dependence
is very strong. Thus, we can image with an adaptive approach that
detects and suppresses the incoherent modes in the data in order to
improve the image.  The high indexed modes are needed for resolution
but they are the first to become incoherent.  Thus, there is a
trade-off between the resolution and robustness of the image, which
leads naturally to an optimization problem solved by the adaptive
approach.  It maximizes a measure of the quality of the image by
weighting optimally the mode amplitudes.

Such mode filtering does not work in waveguides with random media
because there the modes have similar scattering mean free paths.  All
the modes become incoherent at essentially the same propagation
distances and incoherent imaging should be used instead.  There is a
large range interval between the scattering mean free paths of the
modes and the equipartition distance, where incoherent imaging can
succeed. The equipartition distance is the characteristic range scale
beyond which the energy is uniformly distributed between the modes,
independent of the initial state. The waves lose all information about
the source at this distance and imaging becomes impossible.

Incoherent imaging is not useful in waveguides with random boundaries
because the equipartition distance is almost the same as the
scattering mean free paths of the low indexed modes. Once the waves
become incoherent all imaging methods fail.

\section*{Acknowledgements}
We would like to thank Dr Adrien Semin for carrying out the numerical
simulations with Montjoie. The work of L. Borcea was partially supported by the AFSOR Grant
FA9550-12-1-0117, the ONR Grant N00014-12-1-0256 and by the NSF Grants
DMS-0907746, DMS-0934594.  The work of J. Garnier was supported in part
by the ERC Advanced Grant Project MULTIMOD-267184.  The work of C.
Tsogka was partially supported by the European Research Council
Starting Grant Project ADAPTIVES-239959.

\appendix
\section{Numerical simulations of the array data}
\label{sect:ap}

In the numerical simulations the source is supported in a disk of
radius ${\lambda_o}/{10}$, and it emits a pulse
\begin{equation}
f(Bt) = {\rm sinc}\left(B t\right)
\end{equation}
modulated by the carrier signal $\cos(\om_o t)$. The array has $N_R =
39$ receivers located at $\bx_r = (x_r,z_\cA)$, with
$x_r = r  \lambda_o/2$, $r=1,\ldots,39$.

The wave propagation in waveguides with perturbed media is simulated
by solving the wave equation as a first order velocity-pressure
system with the finite element method described in
\cite{becache2000analysis}. It is a second order discretization scheme in space and time,
and in the simulations we used spatial mesh size $h = \la_o/50$  in
cross-range and range,  and time discretization step  determined 
by the CFL condition $\Delta t = {h}/(\sqrt{2} c_{max}),$ with
$c_{max}$ the maximal value of the speed of propagation in the medium.

The wave propagation in waveguides with perturbed pressure release
boundary is simulated by solving the wave equation as a first order velocity-pressure
system  with the code
Montjoie (\url{http://montjoie.gforge.inria.fr/}).  In the simulations
we used $8-$th order
finite elements in space and $4-$th order finite differences in time,
with spatial mesh size $h = \la_{o}/4$ and time discretization step 
$\Delta t = 5\cdot 10^{-6}$\text{s}.

In both cases we use two perfectly matched layers (PML) to the left
and right of the computational domain to model the unbounded waveguide
in $z$.

\end{document}